\documentclass[11pt,draft,reqno]{amsart} 
\usepackage[utf8]{inputenc}
\usepackage{latexsym} 
\usepackage{amssymb} 
\usepackage{amsmath}
\usepackage{enumerate}

\usepackage{hyperref}  

\begin{document}

\newtheorem{theorem}{Theorem}[section]
\newtheorem{problem}{Problem} [section]
\newtheorem{definition}{Definition} [section]
\newtheorem{lemma}{Lemma}[section]
\newtheorem{proposition}{Proposition}[section]
\newtheorem{corollary}{Corollary}[section]
\newtheorem{example}{Example}[section]
\newtheorem{conjecture}{Conjecture} 
\newtheorem{algorithm}{Algorithm} 
\newtheorem{exercise}{Exercise}[section]
\newtheorem{remarkk}{Remark}[subsection]
 
\newcommand{\be}{\begin{equation}} 
\newcommand{\ee}{\end{equation}} 
\newcommand{\bea}{\begin{eqnarray}} 
\newcommand{\eea}{\end{eqnarray}} 

\newcommand{\eeq}{\end{equation}} 

\newcommand{\eeqn}{\end{eqnarray}} 
\newcommand{\beaa}{\begin{eqnarray*}} 
\newcommand{\eeaa}{\end{eqnarray*}}

\newcommand{\lip}{\langle} 
\newcommand{\rip}{\rangle}

\newcommand{\uu}{\underline} 
\newcommand{\oo}{\overline} 
\newcommand{\La}{\Lambda} 
\newcommand{\la}{\lambda} 
\newcommand{\eps}{\varepsilon} 
\newcommand{\om}{\omega} 
\newcommand{\Om}{\Omega} 
\newcommand{\ga}{\gamma} 
\newcommand{\rrr}{{\Bigr )}} 
\newcommand{\qqq}{{\Bigl\|}} 
 
\newcommand{\dint}{\displaystyle\int} 
\newcommand{\dsum}{\displaystyle\sum} 
\newcommand{\dfr}{\displaystyle\frac} 
\newcommand{\bige}{\mbox{\Large\it e}} 
\newcommand{\integers}{{\Bbb Z}} 
\newcommand{\rationals}{{\Bbb Q}} 
\newcommand{\reals}{{\rm I\!R}} 
\newcommand{\realsd}{\reals^d} 
\newcommand{\realsn}{\reals^n} 
\newcommand{\NN}{{\rm I\!N}} 
\newcommand{\DD}{{\rm I\!D}} 
\newcommand{\degree}{{\scriptscriptstyle \circ }} 
\newcommand{\dfn}{\stackrel{\triangle}{=}} 
\def\complex{\mathop{\raise .45ex\hbox{${\bf\scriptstyle{|}}$} 
     \kern -0.40em {\rm \textstyle{C}}}\nolimits} 
\def\hilbert{\mathop{\raise .21ex\hbox{$\bigcirc$}}\kern -1.005em {\rm\textstyle{H}}} 
\newcommand{\RAISE}{{\:\raisebox{.6ex}{$\scriptstyle{>}$}\raisebox{-.3ex} 
           {$\scriptstyle{\!\!\!\!\!<}\:$}}} 
 
\newcommand{\hh}{{\:\raisebox{1.8ex}{$\scriptstyle{\degree}$}\raisebox{.0ex} 
           {$\textstyle{\!\!\!\! H}$}}} 

\newcommand{\OO}{\won} 
\newcommand{\calA}{{\mathcal A}} 
\newcommand{\calB}{{\cal B}} 
\newcommand{\calC}{{\cal C}} 
\newcommand{\calD}{{\cal D}} 
\newcommand{\calE}{{\cal E}} 
\newcommand{\calF}{{\mathcal F}} 
\newcommand{\calG}{{\cal G}} 
\newcommand{\calH}{{\cal H}} 
\newcommand{\calK}{{\cal K}} 
\newcommand{\calL}{{\mathcal L}} 
\newcommand{\calM}{{\mathcal M}} 
\newcommand{\calO}{{\cal O}} 
\newcommand{\calP}{{\cal P}} 
\newcommand{\calU}{{\mathcal U}} 
\newcommand{\calX}{{\cal X}} 
\newcommand{\calXX}{{\cal X\mbox{\raisebox{.3ex}{$\!\!\!\!\!-$}}}} 
\newcommand{\calXXX}{{\cal X\!\!\!\!\!-}} 
\newcommand{\gi}{{\raisebox{.0ex}{$\scriptscriptstyle{\cal X}$} 
\raisebox{.1ex} {$\scriptstyle{\!\!\!\!-}\:$}}} 
\newcommand{\intsim}{\int_0^1\!\!\!\!\!\!\!\!\!\sim} 
\newcommand{\intsimt}{\int_0^t\!\!\!\!\!\!\!\!\!\sim} 
\newcommand{\pp}{{\partial}} 
\newcommand{\al}{{\alpha}} 
\newcommand{\sB}{{\cal B}} 
\newcommand{\sL}{{\cal L}} 
\newcommand{\sF}{{\cal F}} 
\newcommand{\sE}{{\cal E}} 
\newcommand{\sX}{{\cal X}} 
\newcommand{\R}{{\rm I\!R}} 
\renewcommand{\L}{{\rm I\!L}} 
\newcommand{\vp}{\varphi} 
\newcommand{\N}{{\rm I\!N}} 
\def\ooo{\lip} 
\def\ccc{\rip} 
\newcommand{\ot}{\hat\otimes} 
\newcommand{\rP}{{\Bbb P}} 
\newcommand{\bfcdot}{{\mbox{\boldmath$\cdot$}}} 
 
\renewcommand{\varrho}{{\ell}} 
\newcommand{\dett}{{\textstyle{\det_2}}} 
\newcommand{\sign}{{\mbox{\rm sign}}} 
\newcommand{\TE}{{\rm TE}} 
\newcommand{\TA}{{\rm TA}} 
\newcommand{\E}{{\rm E\, }} 
\newcommand{\won}{{\mbox{\bf 1}}} 
\newcommand{\Lebn}{{\rm Leb}_n} 
\newcommand{\Prob}{{\rm Prob\, }} 
\newcommand{\sinc}{{\rm sinc\, }} 
\newcommand{\ctg}{{\rm ctg\, }} 
\newcommand{\loc}{{\rm loc}} 
\newcommand{\trace}{{\, \, \rm trace\, \, }} 
\newcommand{\Dom}{{\rm Dom}} 
\newcommand{\ifff}{\mbox{\ if and only if\ }} 
\newcommand{\nproof}{\noindent {\bf Proof:\ }} 
\newcommand{\nproofYWN}{\noindent {\bf Proof of Theorem~\cite{YWN}:\ }} 
\newcommand{\remark}{\noindent {\bf Remark:\ }} 
\newcommand{\remarks}{\noindent {\bf Remarks:\ }} 
\newcommand{\note}{\noindent {\bf Note:\ }} 
 \newcommand{\examples}{\noindent {\bf Examples:\ }} 
 
\newcommand{\boldx}{{\bf x}} 
\newcommand{\boldX}{{\bf X}} 
\newcommand{\boldy}{{\bf y}} 
\newcommand{\boldR}{{\bf R}} 
\newcommand{\uux}{\uu{x}} 
\newcommand{\uuY}{\uu{Y}} 
 
\newcommand{\limn}{\lim_{n \rightarrow \infty}} 
\newcommand{\limN}{\lim_{N \rightarrow \infty}} 
\newcommand{\limr}{\lim_{r \rightarrow \infty}} 
\newcommand{\limd}{\lim_{\delta \rightarrow \infty}} 
\newcommand{\limM}{\lim_{M \rightarrow \infty}} 
\newcommand{\limsupn}{\limsup_{n \rightarrow \infty}} 
 
\newcommand{\ra}{ \rightarrow } 

 \newcommand{\mlim}{\lim_{m \rightarrow \infty}}  
 \newcommand{\limm}{\lim_{m \rightarrow \infty}}  
 \newcommand{\nlim}{\lim_{n \rightarrow \infty}} 
 
 
 
 
 
 

\newcommand{\one}{\frac{1}{n}\:} 
\newcommand{\half}{\frac{1}{2}\:} 
 
\def\le{\leq} 
\def\ge{\geq} 
\def\lt{<} 
\def\gt{>} 
 
\def\squarebox#1{\hbox to #1{\hfill\vbox to #1{\vfill}}} 
\newcommand{\nqed}{\hspace*{\fill} 
           \vbox{\hrule\hbox{\vrule\squarebox{.667em}\vrule}\hrule}\bigskip}

\title{Average preserving variation processes in view of optimization}


\author{R\'{e}mi Lassalle}

\address{Universit\'{e} Paris 9 (Dauphine), PSL, Place du Mar\'echal De Lattre De Tassigny, 75775 Paris Cedex 16, France}
\email{lassalle@ceremade.dauphine.fr}

\maketitle

\begin{abstract}
In this paper, we investigate specific least action principles for laws of stochastic processes within a  framework which stands on filtrations preserving variations. The associated  Euler-Lagrange conditions, which we obtain, exhibit a deterministic process  in the dynamics aside the canonical martingale term. In particular, taking specific action functionals, extremal processes with respect to those variations encompass specific laws of continuous semi-martingales whose drift characteristic is integrable with independent increments. Then, we relate extremal processes of  classical cost functions, in particular of specific entropy functions, to a class of forward-backward systems of Mckean-Vlasov stochastic differential equations.
\end{abstract}

\section{Introduction}
\label{Section1}

The activity around optimization of functionals over sets of laws of semi$-$martingales, in particular of relative entropy functionals, and around \textit{least action principles}  in stochastic frameworks, covers a wide range of various applications; among many see \cite{1.}, \cite{2.}, \cite{3.}, \cite{4.}, \cite{5.}, \cite{CRESSON}, \cite{6.}, \cite{7.}, \cite{8.}, \cite{9.}, \cite{10.}, \cite{11.}, \cite{MIKAMI}, \cite{CRISAN}, \cite{12.}.   The present paper focusses on \textit{extremal conditions} for functionals depending explicitly on \textit{local characteristics} of laws of \textsc{It\^{o}}'s semi-martingales, under a specific additional constraint. The corresponding least action principle investigated here yields \textsc{Euler-Lagrange} conditions which extend those unconstrained cases investigated in \cite{8.}. For suitable cost maps, it is shown here to yield applications to a specific class of  \textit{forward-backward} systems of coupled   \textsc{Mckean-Vlasov}  stochastic differential equations; the drift terms involve explicit functions of the marginal laws of the underlying processes.   Explicit examples of laws satisfying these conditions are also provided from past dependent \textit{stochastic differential equations}, from a class of processes whose drifts involve some possibly discontinuous \textsc{L\'{e}vy} processes,  and from \textit{partial differential equations}. Within a technical viewpoint, the present paper applies a specific calculus of variations based on \textit{information flows} preserving \textit{transports of measures} (see \cite{8.}, \cite{9.}, \cite{13.}), notably the \textit{intrinsic stochastic derivative} of \cite{13.} : handling the further constraint stated below within this framework requires an accurate control of the induced filtrations. Finally, to provide further insights on this work, we mention here that the results of the present paper will be applied in \cite{Plus1.} to provide extremal conditions to a class of constrained semi-martingale optimization problems (see \cite{12.}), which encompasses \textit{first moment constrained} \textsc{Schr\"{o}dinger} problems; within this specific perspective, as it is expected from the present paper, the extra deterministic process $(A_t^\nu)$ appearing here below may be interpreted as a \textsc{Lagrange} multiplier, which is associated to a further constraint on marginal laws, whose first moment is fixed.

To be accurate, given $d\in \mathbb{N}$, $d\geq 1$, denote by  $W=C([0,1], \mathbb{R}^d)$ the set of continuous $\mathbb{R}^d-$ valued functions on $[0,1]$, which is endowed with the norm of uniform convergence. Subsequently ${\bf M}_1^{\mathbb{S}}(W)$ denotes the subset of laws of continuous \textsc{It\^{o}}'s semi-martingales, such that the evaluation process $({\bf W_t})_{t\in[0,1]}$, which to a given $\omega\in W$ and to a $t\in [0,1]$ associates the value ${\bf W_t(\omega)= \omega(t)}$ of the function $\omega$ at $t$,  is a continuous $\mathbb{R}^d-$ valued process  with the specific form  $${\bf W_t} = {\bf W_0} +  {\bf M_t^\nu} + \int_0^t {\bf v_s^\nu} ds,$$ for all $t\in[0,1]$, $\nu-a.s.$, where ${\bf b^\nu}= \int_0^. {\bf v_s^\nu} ds$ is assumed to be \textit{absolutely continuous} and \textit{adapted} to the $\nu-$\textit{usual augmentation} $(\mathcal{F}_t^\nu)$ of the filtration generated by the evaluation process on $(W, \mathcal{B}_W^\nu, \nu)$, and where $({\bf M_t^\nu})$ is a $(\mathcal{F}_t^{\nu})-$ martingale on the complete probability space $(W, \mathcal{B}_W^\nu, \nu)$, whose covariation $({\bf <M^\nu>_t})_{t\in[0,1]}$ is a $\mathbb{R}^d\otimes\mathbb{R}^d-$ valued $(\mathcal{F}_t^\nu)-$ adapted process of the specific form $${\bf <M^\nu>} = \int_0^. {\bf \alpha_s^\nu} ds, \ \nu-a.s..$$ Furthermore, it is assumed that  $$\int_0^1\left( \|{\bf v_s^\nu}\|_{\mathbb{R}^d}+ \|{\bf \alpha_s^\nu}\|_{\mathbb{R}^d\otimes \mathbb{R}^d}\right) ds <+\infty, \ \nu-a.s..$$ Subsequently,  $({\bf v_t^\nu})$ (respectively $({\bf \alpha_t^\nu})$) will be called the drift (respectively the dispersion) \textit{local characteristic} of $\nu$. We refer to \cite{14.}, \cite{29.}, \cite{30.}, \cite{31.}, \cite{45.}, \cite{PROTTER}, \cite{YORBOOK},  for an introduction on stochastic processes and on martingale theory.  Depending on the context, we may choose a version of $({\bf v_t^\nu})$ which is either \textit{predictable} or \textit{optional} (see \cite{29.}). At any rate, both local characteristics are assumed to be \textit{measurable processes} (see \cite{29.}, section IV). Moreover, ${\bf M}_1^{\mathbb{S}}(W)$ is a subset of the set ${\bf M}_1(W)$ of \textsc{Borel} probability measures on $W$. To provide an insight on the approach unrolled here, it seems to be worth to recall that, for any $q\in C([0,1], \mathbb{R}^d)$ which is $C^1$, we have $\nu= \delta^{Dirac}_q\in {\bf M}_1^{\mathbb{S}}(W)$, the restriction to the \textsc{Borel} $\sigma-$field of the \textsc{Dirac} mass concentrated on $q$; $\delta^{Dirac}_q(A) =1_{A}(q)$, $\forall A\in \mathcal{B}_W$. Subsequently $\lambda_{\mathbb{R}^d}$ (respectively $\lambda|_{[0,1]}$) denotes the \textsc{Lebesgue} measure on $\mathbb{R}^d$, $d\geq 1$ (respectively the restriction of $\lambda_{\mathbb{R}}$ to $[0,1]$).  

 We consider here  functionals $\phi$ whose restriction to ${\bf M}_1^{\mathbb{S}}(W)$ depends explicitly on the \textit{local characteristics} of the laws of semi-martingales,  typically of the form \begin{equation} \phi(\nu) =  \int_{C([0,1],\mathbb{R}^d)\times[0,1]}\mathcal{L}_t({\bf \omega(t)},{\bf v_t^\nu(\omega)},{\bf \alpha_t^\nu(\omega)}) \ \nu\otimes \lambda_{\mathbb{R}}(d\omega,dt), \label{phidef} \end{equation} $\nu \in {\bf M}_1^{\mathbb{S}}(W)$, where $\mathcal{L}$ is measurable and satisfies further technical conditions. Investigations around such functionals are notably motivated by connections to the so-called \textsc{Schr\"{o}dinger}  problem (see \cite{15.}), whose entropic extensions have been related to \textit{optimal transport}  (see \cite{16.}). Recall that, within one of its classical dynamical forms, the \textsc{Schr\"{o}dinger} problem can be seen as a minimization, among \textsc{Borel} probability measures $\nu$ on $C([0,1], \mathbb{R}^d)$, whose initial (respectively final) marginal law  ${\bf W_0}_\star \nu$ (respectively  ${\bf W_1}_\star \nu$) is fixed to be equal to a given \textsc{Borel} probability measure $\nu_0$ (respectively $\nu_1$) on $\mathbb{R}^d$, of the relative entropy $$\mathcal{H}(\nu|\mu_{\mathcal{V}}) = \begin{cases}
        \int_{C([0,1],\mathbb{R}^d)} \ln \frac{d\nu}{d\mu_{\mathcal{V}}}(\omega) \nu(d\omega)& \text{if} \  \nu<< \mu_{\mathcal{V}} \ (i.e. \ absolutely \ continuous) \\
       + \infty & \text{otherwise}
       \end{cases},$$ with respect to a reference measure $\mu_{\mathcal{V}}$; the latter denotes a specific probability measure on $(W,\mathcal{B}_W)$, which is absolutely continuous with respect to 
the \textsc{Wiener} measure $\mu_W$, and whose \textsc{Radon-Nikodym} derivative $\frac{d\mu_{\mathcal{V}}}{d\mu_W}$ is determined by a non-negative differentiable function ${\mathcal{V}}: 
\mathbb{R}^d\to \mathbb{R}$, which satisfies further integrability conditions.  Moreover, within physical models, ${\mathcal{V}}$ is usually identified with a \textit{potential} energy term.
Since at first sight, the relative entropy may seem to be quite distinct from functionals of the form~(\ref{phidef}), as a warm up we briefly recall how \textsc{Schr\"{o}dinger}'s problem actually boils down to the optimization of such functionals.  To start with recall that, following a celebrated sequence of papers of \textsc{R.H. Cameron} and \textsc{W.T. Martin} (see 
\cite{34.}, \cite{35.}), and notably after works of \textsc{G. Maruyama} on \textsc{Markov} processes (\cite{MARUYAMA1}, \cite{MARUYAMA2}) and then of \textsc{I.V. Girsanov} (\cite{GIRSANOV}), classical formulas which determine the representation of the \textsc{Radon-Nikodym} derivatives of absolutely continuous \textsc{Borel} probability measures with respect to the \textsc{Wiener} measure have been obtained, notably by \textsc{T.E. Duncan} and by \textsc{T.Kailath}, where sharp \textsc{It\^{o}}'s stochastic integrals appear as the natural tools to express those fundamental quantities (\cite{DUNCAN}, \cite{KAILATH}, \cite{KZ}, see also \cite{BOOK}). Those formulas yield 
a representation of the relative entropy from \textsc{M. Zakai} et al. (see \cite{BZZ}, formulas (75) and (76)) with respect to the \textsc{Wiener} measure; much recently~\cite{USTUNEL} showed how such results can be used to obtain powerful 
criterions of isomorphisms for specific problems of transformations of the \textsc{Wiener} measure (see \cite
{INV}).  Together with the \textsc{Girsanov} theorem, recall that the latter yields representation formulas as in~\cite{17.} of \textsc{H. F\"{o}llmer} which, together with the \textsc{P.L\'{e}vy} criterion (see \cite{45.}, \cite{PROTTER}),  show that \textsc{Schr\"{o}dinger}'s problem  can be equivalently interpreted as an action functional 
minimization (see \cite{17.}) of the form~(\ref{phidef}) :  \begin{equation} \label{eqvariat} \inf_{\nu \in {\bf M}_1^{\mathbb{S}_B}(W; \nu_0,\nu_1)} \int_{C([0,1],
\mathbb{R}^d)\times [0,1]} \left(\frac{\|{\bf v_t^\nu}({\bf \omega})\|_{\mathbb{R}^d}^2}{2} + \mathcal{V}
(\omega(t)) \right) \nu\otimes\lambda_{\mathbb{R}}({\bf d \omega},dt), \end{equation} $\bf{M_1^{\mathbb{S}_B}}(W;
\nu_0,\nu_1)$ denoting the set $$\left\{\nu\in {\bf M}_1^{\mathbb{S}}(W) :  {\bf W_0}_\star \nu =\nu_0, {\bf W_1}_\star \nu=\nu_1, \alpha^\nu_t(\omega)= I_{\mathbb{R}^d},  \nu
\otimes \lambda|_{[0,1]}-a.e. \right\};$$ the latter functional appearing in~(\ref{eqvariat}) depends explicitly on \textit{local characteristics} of specific laws of \textsc{It\^{o}}'s semi-martingales.

 Therefore, to extend the original \textsc{Schr\"{o}dinger} problem, instead of  using its entropic formulation,  at the inverse one may use the properties of those processes stemming from~(\ref{eqvariat}), and optimize functionals of laws of semi-martingales which depend explicitly on their \textit{local characteristics}; for instance see the \textit{semi-martingale optimal transportation problems} of \cite{12.}, and works of \textsc{T.Mikami} as \cite{MIKAMI}. Recall that a usual \textsc{Schr\"{o}dinger} bridge $(X_t)$ (see \cite{CRUZAMB1}, \cite{17.}, \cite{18.}), which is a stochastic process defined on a \textit{complete stochastic basis} $(\Omega, \mathcal{A}, (\mathcal{A}_t), \mathcal{P})$,  is such that its law $\nu=X_\star\mathcal{P}$ attains the optimum of a problem of the form~(\ref{eqvariat}). Moreover, the marginal density of the absolutely continuous probability measure ${\bf X_t}_\star\mathcal{P}$, with respect to the \textsc{Lebesgue} measure, is of the specific form $$\rho(t,{\bf x})= \theta(t, {\bf x}) \theta^\star(t, {\bf x}), \ \lambda_{\mathbb{R}^d}-a.e.,$$ i.e. $$\mathcal{P}(X_t \in A)= \int_{A} \theta(t, {\bf x}) \theta^\star(t, {\bf x})\lambda_{\mathbb{R}^d}(d{\bf x}),$$ for any \textsc{Borel} set $A\in \mathcal{B}_{\mathbb{R}^d}$, and $t\in(0,1)$. Furthermore, $\theta$ satisfies a \textit{partial differential equation}, which is usually interpreted as an \textit{Euclidean} version of the \textsc{Schr\"{o}dinger} equation (see \cite{CRUZAMB1}), since it is usually of the form 
$$\begin{cases}-\frac{\partial \theta^\star}{\partial t} = \mathcal{H} \theta^\star  \\
       \theta^\star(0, {\bf x}) = f({\bf x}), \end{cases}$$ where $$\mathcal{H}=- \frac{\Delta}{2} + \mathcal{V},$$  where $\mathcal{V}$ is a nice function which satisfies further technical conditions (see \cite{CRUZAMB1}),  and where $f\in \Dom(H)$. Similarly, $\theta$ satisfies an adjoint equation (see \cite{CRUZAMB1}, \cite{18.}). Using those properties, the specific value of the \textsc{Radon-Nikodym} derivatives of the law of those usual \textsc{Schr\"{o}dinger} bridges $(X_t)$ with respect to the \textsc{Wiener} measure can be computed explicitly from the \textsc{Clark-Bismut-Ocone} formula of the \textsc{Malliavin} calculus (see \cite{33.}). From this explicit value, works in the line of \cite{18.} (in particular see \cite{CRUZAMB1})  have shown that the associated stochastic processes $(X_t)$ satisfied strong extremal conditions on the fixed stochastic basis of the \textsc{Wiener} space. This yields a critical variational property of those strong extremal points among a class of \textsc{Bernstein}'s (reciprocal) processes realized on a fixed complete stochastic basis, in this case the \textsc{Wiener} space (see \cite{CRUZAMB1}).  Note that to be applied, the previous method requires to first obtain precise information on explicit properties of the law of the process which attains~(\ref{eqvariat}). However, consider now the case of more general optimization problems than~(\ref{eqvariat}), still with fixed initial and final marginals, but allowing various specific constraints on the dispersion \textit{local characteristic}, and where many general cost functionals are now considered, typically of the form~(\ref{phidef}); it still depends explicitly on the \textit{local characteristics} of laws of semi-martingales, and in some cases  such functionals may also be allowed to depend on further features of the law, as the density of intermediate marginals. In these latter cases, contrary to what happens in the case~(\ref{eqvariat}) under suitable hypothesis on $\mathcal{V}$,  the laws of the optima are not necessarily known explicitly from calculus, neither the explicit value of their drifts. Moreover, for certain values of the dispersion \textit{local characteristic}, and notably in the singular case, not encompassed by~(\ref{eqvariat}), where the dispersion vanishes, the law of the optima are not necessarily absolutely continuous with respect to the \textsc{Wiener} measure, and may not allow to apply the sharp \textit{stochastic calculus of variations} of \textsc{P.Malliavin} (see \cite{33.}) to perform computations; indeed the involved variations may not satisfy the required \textit{quasi-invariance} properties.  Due to the specificities of those problems, to obtain precise extremal conditions together with an efficient variational principle, which applies indistinctly to probability measures concentrated on paths which satisfy the classical \textsc{Euler-Lagrange} condition, and to extensions of the \textsc{Schr\"{o}dinger} problem, it suggests to work in a framework (see~\cite{8.}) which allows to perturbate the stochastic basis itself rather than the process, and to perform an extension of the mathematical version of the so-called \textsc{Hamilton}'s least action principle, originated from physics (see \cite{21.}, \cite{22.}, \cite{CARTAN} p.293, \cite{24.}, \cite{23.}); recall that, in the classical case, the latter doesn't require to know the explicit value of the optima to deduce the extremal conditions.  This lead \cite{8.} to consider extremal conditions of the form \begin{equation} \bf{\partial}_{\bf v}
\mathcal{L}_t({\bf W_t},{\bf v_t^\nu},{\bf \alpha_t^\nu}) - \int_0^t {\bf \partial}_{\bf x}\mathcal{L}_s({\bf W_s},{\bf v_s^\nu},{\bf \alpha_s^\nu}) ds = {\bf N_t^\nu}, \label{ISELAPor}\end{equation} up to some negligible sets,  where $({\bf N_t^\nu})$ is a \textit{c\`{a}d-l\`{a}g}   $(\mathcal{F}_t^\nu)-$martingale; using the conventional terminology of stochastic analysis (see \cite{29.}, \cite{30.}), we call \textit{c\`{a}d-l\`{a}g} the right-continuous functions with left limits and similarly for stochastic processes. Note that conditions as~(\ref{ISELAPor}) require that a certain process which depends explicitly on the \textit{local characteristics} of the law is a \textit{martingale}. This explicit dependency allows to make a distinction between, on the one hand the not necessarily convex problem to determine laws which satisfy conditions as~(\ref{ISELAPor}), and on the other hand usual martingale problems on a canonical space (see \cite{JACODS}). Relating~(\ref{ISELAPor}) to critical conditions of functionals  as ~(\ref{phidef}), through a \textit{least action principle} under general hypotheses, requires to extend to this specific framework the classical \textit{calculus of variations}, and to perform specific variations of functionals on ${\bf M}_1(W)$; for an introduction to the basics of \textit{calculus of variations}, we refer here to the classical lecture notes of \textsc{H.Cartan} (\cite{CARTAN}). Since the functionals of interest depend explicitly on \textit{local characteristics} of laws of semi-martingales, this emphasizes phenomena of \textit{filtration theory} which are specific to the existence of a \textit{noise} in stochastic frameworks. Those aspects are related to the \textit{variation of local characteristics along specific transports of measure} (see \cite{9.}). This provides a certain geometric point of view on the \textit{innovation conjecture}  of \textit{filtering} (see \cite{25.}, \cite{27.},  \cite{26.}, \cite{ABSAPI}, \cite{MEYERF}, \cite{USTUNEL}, \cite{USTUNEL2}) which is closely related to the mysterious  \textsc{B.Tsirelson}'s counter-example (see \cite{28.}), and to the counter-example of \cite{DUBINS};  among many insights on \cite{DUBINS}, see also \cite{BEG} and the references therein, and recall that it stands on the \textsc{A.Vershik} theory of filtration (see also \cite{EMSC}). The same aspects of regularity of local characteristics due to those phenomena also occur for variations preserving initial and final marginals, which may be seen though localization arguments as those of  \cite{LASASU}, \cite{ABSAPI}, \cite{Loc}.   The construction which is used in this paper actually stands on variations which satisfy a stronger constraint as the innovation conjecture, and is therefore justified by difficulties which occur from the fact that some semi-martingales can not be determined from their innovation process (see \cite{MEYERF}, \cite{Tsirelinnov}).   From this, and since for reasons explained above, the \textsc{Malliavin} calculus is not sufficient to achieve this purpose,  it seems necessary to complete those tools of \textit{stochastic analysis} (\cite{33.}) by another calculus on the set of laws of stochastic processes, to incorporate those sharp specificities inherent to \textit{filtration theory}, while keeping a compact notation. Such a construction, which relies on \textit{information flows} preserving transports of measure has been proposed in  \cite{13.}, which will be applied throughout this paper. As it is developed on the canonical space endowed with a filtration, the corresponding derivative is not conditioned on the choice of a stochastic basis, where a specific model would be considered; we call it the \textit{intrinsic stochastic derivative}.  From those latter specific features, this must be distinguished from works in the line of \cite{20.}. Within this framework of \cite{13.}, the present paper investigates  \textsc{Euler-Lagrange} conditions which, outside some negligible sets, are of the form \begin{equation} \bf{\partial}_{\bf v}
\mathcal{L}_t({\bf W_t},{\bf v_t^\nu},{\bf \alpha_t^\nu}) - \int_0^t {\bf \partial}_{\bf x}\mathcal{L}_s({\bf W_s},{\bf v_s^\nu},{\bf \alpha_s^\nu}) ds = {\bf N_t^\nu} +{\bf A_t^\nu} , \label{ISELAP}\end{equation} where $({\bf N_t^\nu})_{t\in[0,1)}$ is a \textit{c\`{a}d-l\`{a}g}   $
(\mathcal{F}_t^\nu)-$martingale, and where $({\bf A_t^\nu})_{t\in[0,1)}$ is a deterministic process, together with the associated least action principle. Taking specific cost functions, this yields information on some laws of semi-martingales whose drift characteristic is an integrable process with independent increments. By taking $A_t^\nu=0$, for all $t\in[0,1)$, note that~(\ref{ISELAP}) extends~(\ref{ISELAPor}) of \cite{13.}.

The structure of this paper is the following. In Section~\ref{Section2}, we introduce the notation, used in the whole paper. In Section~\ref{Section3}, we recall accurately the specific tools, notably those of  \cite{13.}, which are used subsequently to obtain a compact notation, taking into account the circumstances identified above, which we encounter within the specific context of those filtered probability spaces; several counter-examples are also recalled. In Section~\ref{Section4}, we introduce average preserving variation processes (Definition~\ref{averagepreservingdef}), and investigate those of their properties which we use subsequently. It enables us to establish in Section~\ref{Section5} our main result Theorem~\ref{aplapT}, an average preserving least action principle, which characterizes those laws of semi-martingales which satisfy~(\ref{ISELAP}), as extremal points of functionals of the form~(\ref{phidef}), with respect to a set of average preserving variations. In Section~\ref{Section6},  for classical cost functions, Proposition~\ref{PropMVFB} relates critical points of the average preserving least action principle, to a specific class of \textit{forward-backwards} systems (for instance see \cite{32.}) of \textsc{Mckean-Vlasov} stochastic differential equations. Finally, three explicit examples are provided.

\section{Notation}
\label{Section2} 
Given $d\in \mathbb{N}$, with $d\geq 1$, $W=C([0,1], \mathbb{R}^d)$ denotes the set of continuous $\mathbb{R}^d-$valued functions on $[0,1]$; $W$ stands for \textsc{Wiener}, whose probability space provides here our main paradigm. It is endowed with the norm $\|.\|_W$ of uniform convergence, which turns it into a separable \textsc{Banach} space, whose related \textsc{Borel} sigma-field is denoted by $\mathcal{B}_W$. The set of \textsc{Borel} probability measures on $W$ is denoted by ${\bf M}_1(W)$. Given $\eta\in {\bf M}_1(W)$, $\mathcal{B}_W^\eta$ denotes the $\eta-$completion of the \textsc{Borel} sigma-field. 

If $(\Omega, \mathcal{A}, \mathcal{P})$ is a probability space, and $X :\Omega \to \mathbb{R}$ is an $\mathcal{A}/\mathcal{B}_{\mathbb{R}}-$ measurable function, we use the standard notation $$\mathbb{E}_{\mathcal{P}}\left[X\right]= \int_\Omega X(\omega) \mathcal{P}(d\omega)$$ of the mathematical expectation, whenever $X\geq 0$, $\mathcal{P}-a.s.$ in which case $\mathbb{E}_{\mathcal{P}}\left[X\right]\in \mathbb{R}_+\cup\{+\infty\}$, or when $\mathbb{E}_{\mathcal{P}}[|X|]<+\infty$, in which case $\mathbb{E}_{\mathcal{P}}[X]\in \mathbb{R}$ and $X$ is said to be $\mathcal{P}-$ integrable; $\mathcal{B}_{\mathbb{R}}$ denotes the usual \textsc{Borel} $\sigma-$field on $\mathbb{R}$. 

The classical \textsc{Cameron-Martin} space $H^1$ is defined by $$H^1= \left\{ {\bf h}: [0,1] \to \mathbb{R}^d : \ {\bf h}=\int_0^. {\bf \dot{h}_s} ds \ , \ \int_0^1\|{\bf \dot{h}_s}\|^2_{\mathbb{R}^d}ds <+\infty \ \right\}.$$ Recall that it is turned into a \textsc{Hilbert} space with product $$ <{\bf h},{\bf k}>_{H^1}  =  \int_0^1 <{\bf \dot{h}_s},{\bf \dot{k}_s}>_{\mathbb{R}^d} ds ,$$ for all ${\bf h},{\bf k}\in H^1$; $\|{\bf h}\|_{H^1} = \sqrt{<{\bf h},{\bf h}> _{H^1}}$.  This space plays a key role in stochastic analysis (see \cite{34.}, \cite{35.}, \cite{33.}). However in view of establishing least action principles, its vector subspace  $$H^{\bf 1}_{0,0} = \left\{ {\bf h}\in H^1 : {\bf h_0}={\bf h_1}={\bf 0}_{\mathbb{R}^d} \right\},$$ will play a major role.

\section{Recall on the intrinsic stochastic derivative on ${\bf M}_1(W)$} 
\label{Section3}
\subsection{Information flows preserving transports of measure}

Let $\eta \in {\bf M}_1(W)$, we denote by $M_{\eta}((W,\mathcal{B}_W^\eta),(W,\mathcal{B}_W)$) the set which is obtained by identifying the  $\mathcal{B}_W^\eta / \mathcal{B}_W-$ measurable functions ${\bf f} : W\to W$, which coincide outside an $\eta-$negligible set. Given ${\bf U}\in M_{\eta}((W,\mathcal{B}_W^\eta),(W,$ $\mathcal{B}_W)$), the filtration $(\mathcal{G}_t^{\bf U})$ generated by ${\bf U}$ is the $\eta-
$usual augmentation of the filtration $(\sigma({\bf f_s}, s\leq t))_{t\in[0,1]}$, for some (and then all) $\mathcal{B}_W^\eta / \mathcal{B}_W-$ measurable ${\bf f}: W\to W$ whose $\eta-$ equivalence class is ${\bf U}$, where ${\bf f_s}= {\bf W_s}\circ {\bf f}$, $s\in[0,1]$; $({\bf W_s})$ denotes the evaluation process, recall that ${\bf W_t}({\bf \omega})= {\bf \omega}(t)$, for all $t\in[0,1]$, ${\bf \omega} \in W$. We denote by $(\mathcal{F}_t^\eta)$ the $\eta-$usual augmentation of the filtration generated by the evaluation process on $(W, \mathcal{B}_W^\eta, \eta)$. A  ${\bf U}\in M_{\eta}((W,\mathcal{B}_W^\eta),(W,$ $\mathcal{B}_W)$) is said to be $(\mathcal{F}_t^\eta)-$adapted if $$(\mathcal{G}_t^{\bf U}) \subset (\mathcal{F}_t^\eta) ;$$ since $(\mathcal{F}_t^\eta)$ satisfies the usual conditions, it is equivalent to ${\bf f_s}$ is $\mathcal{F}_s^\eta-$measurable, for all $s\in[0,1]$,  for some (and then all)  $\mathcal{B}_W^\eta / \mathcal{B}_W-$ measurable ${\bf f} : W\to W$, whose $\eta-$equivalence class is ${\bf U}$.

We say that ${\bf U}\in M_{\eta}((W,\mathcal{B}_W^\eta),(W,$ $\mathcal{B}_W))$  is an isomorphism of filtered probability space,  if it is $(\mathcal{F}_t^\eta)-$adapted, and if there exists a $(\mathcal{F}_t^\nu)-$adapted
${\bf \widetilde{U}}\in M_\nu((W,\mathcal{B}_W^\nu)$, $(W,\mathcal{B}_W)$), 
where $\nu={\bf U}_\star \eta$, which is such that $${\bf \widetilde{U}}\circ {\bf U} = {\bf I}_W,  \ \eta-a.s.,$$ and $$ {\bf U}\circ {\bf \widetilde{U}} = {\bf I}_W,  \ \nu-a.s.,$$  ${\bf I}_W: \omega\in W \to \omega \in W$ denoting the identity map on $W$, and ${\bf U}_\star \eta$ denoting the direct image (push-forward of measure) of $\eta$ by ${\bf U}$; we call ${\bf \widetilde{U}}$ the inverse of ${\bf U}$. Such isomorphisms are well known fundamental mathematical objects in stochastic analysis (see \cite{33.}, \cite{Tsireltriple}). Recall that isomorphisms of filtered probability spaces are exactly information flows preserving maps, in the acceptation that ${\bf U}\in M_{\eta}
((W,\mathcal{B}_W^\eta),(W,\mathcal{B}_W))$ is an isomorphism of filtered probability space, if and only if,  $$(\mathcal{G}_t^{\bf U}) = (\mathcal{F}_t^\eta) ;$$ see \cite{13.}. Here we use the term information flows within the terminology of \cite{14.} p.39, which interprets a \textit{filtration} as an \textit{information flow}.

 Subsequently, since $
\sigma({\bf W_0})^\eta$ (the $\eta-$completion of the $\sigma-$field $\sigma({\bf W_0})$) does not necessarily coincide with $\mathcal{F}_0^\eta$, to obtain an efficient set of variation processes,  it is useful to introduce 
the subset $\mathcal{I}_f^0(\eta)$ of the ${\bf U}\in M_{\eta}((W,\mathcal{B}_W^\eta),(W,$ $\mathcal{B}_W))$ which are isomorphisms of filtered spaces, and further satisfy $$\sigma
({\bf W_0})^\eta = \sigma({\bf U_0})^\eta.$$ It is enlightening to interpret $\mathcal{I}_f^0(\eta)$ as the set of  information flows preserving maps on $(W,\mathcal{B}_W^\eta,\eta)$ which also 
preserve the initial information.

Recall that transport plans (see \cite{36.}) notably enable relaxations of the \textsc{Monge} problem~\cite{37.} (see \cite{38.}, \cite{36.}). Let $E$ (respectively $S$) be Polish spaces endowed with filtrations $(\mathcal{B}_{t,E})_{t\in I}$ (respectively $(\mathcal{B}_{t,S})_{t\in I}$) of their \textsc{Borel} sigma-fields, labeled by a same totally ordered set $I$. Given $\eta \in {\bf M}_1(E)$, $\nu\in {\bf M}_1(S)$, the set $\Pi(\eta,\nu)$ of transport plans  of $\eta$ to $\nu$, is the set of $\gamma \in {\bf M}_1(E\times S)$ whose first (respectively second) marginal is $\eta$ (respectively $\nu$). Denote by  $Q_\gamma$ the conditional probability kernel (see \cite{29.}, \cite{DM3}), such that $\gamma = \int_E \eta(dx) \delta^{Dirac}_x\otimes Q_\gamma^x$. Within those hypothesis, recall that any such $\gamma$ generates a filtration $(\mathcal{G}_t(\gamma))$ on the measurable space $(E, \mathcal{B}_E^\eta)$, such that for all $t\in I$, $\mathcal{G}_t(\gamma)$ is the $\eta-$completion of the smallest sigma-field such that for all $B\in \mathcal{B}_{t,S}$ of $\nu-$continuity (i.e. $\nu(\partial B)= 0$), the map $\phi_B : x\in E \to Q_\gamma^x(B) \in[0,1]$ is measurable. 
\begin{definition}\label{def31ISOPLAN}
 Let $\eta \in {\bf M}_1(E)$, $\nu\in {\bf M}_1(S)$. We say that $\gamma \in \Pi(\eta,\nu)$ is an information flows preserving transport plans from $\eta$ to $\nu$, if it further satisfies $$\mathcal{G}_t(\gamma) = \mathcal{B}_{t,E}^\eta, \ \forall t\in I,$$ $(\mathcal{G}_t(\gamma))_{t\in I}$ denoting the filtration generated by $\gamma$ on $(E, \mathcal{B}_E^\eta)$.
\end{definition}

Set $\mathcal{B}_{t,W}^0= \sigma(W_s, s\leq t)$, take $E= S=W$, $I=[0,1]$, $\eta \in {\bf M}_1(W)$, and  take
$\mathcal{B}_{t,E}=\mathcal{B}_{t+,W}^0$ and $\mathcal{B}_{t,S} = \mathcal{B}_{t,W}^0$, for any $t\in[0,1]$. Then,  isomorphisms of filtered probability spaces ${\bf U}$ on $(W,\mathcal{B}_E^\eta, \eta)$  such that $\nu= {\bf U}_\star\eta$ are identified with deterministic transport plans of this kind.

\begin{remark}
Although in several cases, for $t\in I$,  $\mathcal{G}_t(\gamma)$ coincides with $\sigma(\phi_B : B \in \mathcal{B}_t(S))^\eta$, it is not necessarily always the case. For the reader's 
convenience, we recall here a simple elementary counter-example based on not finite, but countable, spaces. Let $E$ and $S$ be the \textsc{Borel} subset of $\mathbb{R}$ 
defined by $$E=S =\{0,1\} \bigcup \left\{\frac{1}{n+1} : n\in \mathbb{N} \right\},$$ and endow each of those spaces with the Polish topology induced by the restriction of the usual 
Euclidean norm to those subsets; take $I=\{0,1\}$ and set $\mathcal{B}_{0,E} = \mathcal{B}_{0,S}=\{\emptyset, S\}$, $\mathcal{B}_{1,E}= \mathcal{B}_{1,S}= \{\emptyset, \{0\}, S
\setminus\{0\}, S\} \subset \mathcal{B}_S$. From a \textsc{Poisson} distribution on the set of natural numbers, we may obtain a \textsc{Borel} probability measure $\eta\in {\bf M}_1(E)$ 
on the infinite countable space $E$, which further satisfies $\eta(\{x\})>0$, $\forall x\in E$. Let $(I_E, T) : x\in E \to (x, T(x)) \in E\times S$, where ${T} : E \to S$ is the \textsc{Borel} measurable function defined by $T(x)=1_{E\setminus \{0\}}(x) \in S$, $\forall x\in E$.  Since $(I_E, T)$ is \textsc{Borel} measurable, the pushforward $\gamma =(I_E,T)_\star \eta$ of $\eta$ by $(I_E,T)$ is a well defined element of ${\bf M}_1(E\times S)$, which is the deterministic transport plan of $\eta$ induced by $T$. From the definitions, we have $$\gamma(A\times B) = \mathbb{E}_{\eta}[1_A Q(B)],$$ $\forall A\in \mathcal{B}_E$, $ \forall B\in \mathcal{B}_S$, where $Q : E \to {\bf M}_1(S)$ is given by $$Q^x = \delta^{Dirac}_{T(x)}, \ \forall x\in E, $$ $\delta^{Dirac}_y
$ denoting the \textsc{Dirac} mass centered on $y$, $\forall y\in S$; we adopt the same conventional terminology as section 5.4.1 of \cite{AirMal}. On the other hand, from the hypothesis, the 
definitions yield $\gamma \in \Pi(\eta, \nu)$, where $\nu=T_{\star}\eta \in M_1(S)$ is necessarily of the form $$\nu = p\delta^{Dirac}_0 + (1-p) \delta^{Dirac}_1,$$ for some $p\in(0,1)$.  
Then, it is an easy task to check that $\mathcal{G}_1(\gamma)$ is strictly smaller than $\sigma(\phi_B : B\in \mathcal{B}_{1,S})^\eta=\mathcal{B}_{1,E}$, for $\{0\}\in \mathcal{B}_
{1,E}$ but  $\{0\} \notin \mathcal{G}_1(\gamma)$.  (End of Remark.) \end{remark}

\subsection{Variation processes.}

\textit{Variation processes} have been defined in \cite{13.}. Recall that given $\eta \in {\bf M}_1(W)$, $L^2_a(\eta,H^1)$ is the subset of the ${\bf h} \in M_{\eta}((W,\mathcal{B}_W^\eta),(W,\mathcal{B}_W))$, such that $$\mathbb{E}_\eta[\|{\bf h}\|_{H^1}^2]<+\infty,$$ which are further assumed to be $(\mathcal{F}_t^\eta)-$adapted.  It is an \textsc{Hilbert} space with product $$<{\bf h},{\bf k}>_{L^2_a(\eta,H^1)}= \mathbb{E}_\eta\left[<{\bf h},{\bf k}>_{H^1}\right] = \mathbb{E}_\eta\left[\int_0^1<{\bf \dot{h}_s},{\bf \dot{k}_s}>_{\mathbb{R}^d}ds \right],$$ for all ${\bf h},{\bf k}\in L^2_a(\eta,{H^1})$. Given ${\bf h}\in L^2_a(\eta,{H^1})$, we set ${\bf \tau_h} = {\bf I}_W+ {\bf h}$, ${\bf I}_W$ still denoting the identity map on $W$.

Further recall that the set of \textit{variation processes} $V_\eta$  at $\eta$, which may be interpreted as the set of perturbations preserving information flows, is the subset  of the ${\bf h} \in L^2_a(\eta, {H^1})$ which are ruled by the following \textit{absence of information loss} principle : $${\bf U} \in \mathcal{I}_f^0(\eta) \implies {\bf U}+ {\bf h} \in \mathcal{I}_f^0(\eta), \ \forall {\bf U} \in M_{\eta}((W,\mathcal{B}_W^\eta),(W, \mathcal{B}_W)).$$ 
\begin{remark}
For $\eta \in M_1(W)$, ${\bf h}\in L^2_a(\eta,{H^1})$, and ${\bf U} \in \mathcal{I}_f^0(\eta)$, we have $$(\mathcal{G}_t^{{\bf U}+ {\bf h}}) \subset (\mathcal{F}_t^\eta).$$ However, the previous inclusions may be strict, in which case ${\bf U}+ {\bf h} \notin \mathcal{I}_f^0(\eta)$. For instance, with $d=1$,  let $\eta$ be the law of weak solutions to the so-called \textsc{Tsirelson}'s equations (see notably \cite{YORBOOK}, Chapter IX, p. 392, or \cite{45.}), let ${\bf U}={\bf I_W}$ (identity map on $W$), and set ${\bf h}= -\int_0^. v_t(\omega) dt\in L^2_a(\eta,H^1)$, $v$ denoting \textsc{Tsirelson}'s drift. From the \textsc{Girsanov} theorem (\cite{GIRSANOV}), it is known that $\left({\bf U}_t+ {\bf h}_t\right)_{t\in[0,1]}$ is an $(\mathcal{F}_t^\eta)-$ Brownian motion. However, the filtration generated by ${\bf U}+ {\bf h}$ is strictly smaller than $(\mathcal{F}_t^\eta)$, since \textsc{Tsirelson}'s equation has no strong solutions (\cite{28.}, \cite{45.}). As a consequence, ${\bf U}+ {\bf h}\notin  \mathcal{I}_f^0(\eta)$. Stronger counter-examples may also be obtained from the so-called second \textsc{Tsirelson}'s counter-example (see \cite{DUBINS}), and with localization arguments. (End of Remark.)
\end{remark}

 It follows from the definition that $V_\eta$ is a linear subspace of $L^2_a(\eta,H^1)$ (see \cite{13.}). In view of applications to optimization, and to least action principles, the following sets are useful : $$V_\eta^\infty= \left\{ {\bf h}\in V_\eta : \exists \ C>0 : \|{\bf h}\|_W \leq C, \ \eta-a.s. \ \right\},$$  $$V_\eta^{0,\infty} = V_\eta^\infty \cap L^2_a(\eta, H^{\bf 1}_{0,0}),$$ i.e.  ${\bf h}\in V_\eta^\infty$ is an element of $V_\eta^{0,\infty}$, if and only if, ${\bf h_0}={\bf h_1}={\bf 0}_{\mathbb{R}^d}, \ \eta-a.s.$;  recall that from Proposition 2.4  (respectively from Lemma 2.1) of \cite{13.},  although, depending on $\eta$,  $V_\eta$ and $V_\eta^\infty$ (respectively $V_\eta^{0,\infty}$) are not necessarily closed, they are however dense in $L^2_a(\eta,H^1)$ (respectively in $L^2_a(\eta,H^{\bf 1}_{0,0})$), $\forall \eta \in {\bf M}_1(W)$.

\begin{remark}
For the reader's convenience, recall that for $\eta \in {\bf M}_1^{\mathbb{S}}(W)$, ${\bf h}\in V_\eta$, and $\epsilon \in \mathbb{R}$, we have not only $\eta^{\epsilon {\bf h}} = {{\bf \tau}_{\epsilon {\bf h}}}_\star\eta \in {\bf M}_1^{\mathbb{S}}(W)$, as it is expected from the so-called \textsc{Stricker}'s theorem (see Theorem 4 of \cite{PROTTER}) on \textit{filtrations shrinkage}, but also ${\bf b^{\eta^{\epsilon {\bf h}}}}\circ {\bf \tau}_{\epsilon {\bf h}} = {\bf b^\eta} + \epsilon {\bf h}$, $\eta-a.s.$, $\forall \epsilon \in \mathbb{R}$. This latter property justifies, in part, the definition of \cite{13.} recalled in the next subsection. (End of Remark.) 
\end{remark}

\subsection{The intrinsic stochastic derivative.}

The \textit{intrinsic stochastic derivative} has been defined in \cite{13.}, in view of differentiating functionals with finite values on subsets of ${\bf M}_1^{\mathbb{S}}(W)$, which depend explicitly on \textit{local characteristics}. It is motived by the variation of \textit{local characteristics} along adapted transports of measure on ${\bf M}_1^{\mathbb{S}}(W)$, and permits to obtain directly compact statements.

 We recall part of its definition, which is necessary to state the least action principle with average preserving variations of Section~5. Given a function $$\phi : \eta \in {\bf M}_1(W) \to \phi(\eta) \in \mathbb{R}\cup\{+\infty\},$$ and $\eta \in {\bf M}_1(W)$ such that $\phi(\eta)<+\infty$, $\phi$ is said to be $L^2_a(\eta, H^{\bf 1}_{0,0})-$differentiable at $\eta$ if for all ${\bf k}\in V_\eta^{0,\infty}$,  $
 \frac{d}{d\epsilon}\phi(\eta^{\epsilon {\bf k}})\big|_{\epsilon =0}$ exists, where $$\eta^{\epsilon {\bf k}}= ({\bf I}_W+\epsilon {\bf k})_\star \eta,$$ for all $\epsilon \in \mathbb{R}$, and if there exists ${\bf \xi} \in L^2_a(\eta, H^{\bf 1}_{0,0})$ such that $$
 \frac{d}{d\epsilon}\phi(\eta^{\epsilon {\bf k}})\big|_{\epsilon=0} = \int_W <{\bf \xi}, {\bf k}>_{H^1} d\eta,$$ for all ${\bf k}\in V_\eta^{0,\infty}$. In this case we define $$\delta \phi_\eta : {\bf k} \in 
 L^2_a(\eta,H^1_{0,0}) \to \int_W <{\bf \xi}, {\bf k}>_{H^1} d\eta \in \mathbb{R}.$$ Note that, by definition we have $\delta \phi_\eta[{\bf k}]= \frac{d}{d\epsilon}\phi(\eta^{\epsilon {\bf k}})\big|_{\epsilon =0},$ for all ${\bf k}\in V_\eta^{0,\infty}$, which motivates an efficient definition of \textit{variation processes}.

\section{Average preserving variations processes.} 
\label{Section4}
\begin{definition}\label{averagepreservingdef}
For all $\nu \in {\bf M}_1(W)$, we define the set of average preserving variation processes to be the set $$A_\nu^{0,\infty}= \left\{ {\bf h} \in V^{0,\infty}_\nu : \int_{W} {\bf h} \ d\nu= {\bf 0}_{H^1}\right\}, $$ where $\int_W {\bf h} \ d\nu$ is a \textsc{Bochner} integral (see \cite{39.}).
\end{definition}
\begin{remark}
Since $V_\nu^\infty \subset L^2_a(\nu,H^1)$, $\int_W {\bf h} \ d\nu$ is a well defined element of the separable \textsc{Hilbert} space $H^1$. (End of Remark.)
\end{remark}

\begin{proposition} Let $\nu\in {\bf M}_1(W)$, and $$ j :  {\bf h} \in L^2_a(\nu, H^1) \to {\bf h}-\int_W {\bf h} \ d\nu \in L^2_a(\nu,H^1), $$ then we have $A_\nu^{0,\infty}= j(V^{0,\infty}_\nu).$ In particular $A_\nu^{0,\infty}$ is a vector space.  \end{proposition}
{\bf Proof :}
The inclusion $A_\nu^{0,\infty} \subset j(V^{0,\infty}_\nu)$ follows from the definitions. Conversely, let ${\bf k}\in j(V^{0,\infty}_\nu)$ and ${\bf h}\in V^{0,\infty}_\nu $ be such that ${\bf k}= {\bf h}- 
\int_W {\bf h} \ d\nu$. The average of ${\bf k}$ vanishes, and since ${\bf h}$ is essentially bounded, the same holds with ${\bf k}$. Moreover since ${\bf h_0}={\bf h_1}={\bf 0}_{\mathbb{R}^d},$ $\nu-a.s.$, we obtain  ${\bf k_0}={\bf k_1}={\bf 0}_{\mathbb{R}^d}, \ \nu-a.s.,$ so that it is enough to prove that ${\bf k}\in V_\nu$. Let ${\bf U}\in \mathcal{I}_f^0(\nu)$ be an isomorphism of filtered 
probability space on $(W, \mathcal{B}_W^\nu,\nu)$. We set ${\bf U}^{\bf k}= {\bf U}+ {\bf k}.$ First notice that since ${\bf k_0}={\bf 0}_{\mathbb{R}^d}, \  \nu-a.s.$, we have $$\sigma({\bf U_0^k})^\nu= \sigma({\bf U_0})^\nu = \sigma({\bf W_0})^\nu.$$ Henceforth, to avoid a heavy notation in this proof, we set $$\mathbb{E}_\nu\left[{\bf h}\right] = \int_W {\bf h} d\nu.$$
By definition we have $ {\bf U}^k=  {\bf T} +{\bf h},$ where ${\bf T} =  \tau_{-\mathbb{E}_\nu\left[{\bf h}\right]}\circ {\bf U}, \ \nu-a.s.,$ and where $$  \tau_{-\mathbb{E}_\nu\left[{\bf h}\right] } : {\bf \omega} \in W \to {\bf \omega}-\mathbb{E}_\nu
[{\bf h}] \in W.$$ Note that ${\bf T}$ is well defined, that it is $(\mathcal{F}_t^\nu)-$adapted, and that it can be checked to satisfy $(\mathcal{G}_t^{{\bf T}})=(\mathcal{G}_t^{\bf U})$. Therefore, we obtain that ${\bf T} \in 
\mathcal{I}_f^0(\nu)$. Since ${\bf h}\in V_\nu$, ${\bf T}\in \mathcal{I}_f^0(\nu)$ implies ${\bf U}^k\in \mathcal{I}_f^0(\nu)$. Thus, 
for all ${\bf U}\in \mathcal{I}_f^0(\nu)$, ${\bf U}^k = {\bf U}+{\bf k} \in \mathcal{I}_f^0(\nu).$ This shows that ${\bf k}\in V_\nu$, and thus, that ${\bf k} \in V^{0,\infty}_\nu$ with $\int_{W} {\bf k} d\nu={\bf 0}_{H^1}$. Whence $ j(V^
{0,\infty}_\nu)\subset A_\nu^{0,\infty}.$  \qed

\begin{proposition} For any $\nu\in {\bf M}_1(W)$, the closure $cl(A_\nu^{0,\infty})$ of $A_\nu^{0,\infty}$ in $L^2_a(\nu,H^1)$ satisfies $$cl(A_\nu^{0,\infty}) = \left\{ {\bf h}\in L^2_a(\nu, H^{\bf 1}_{0,0}) :  \int_W {\bf h} \ d\nu ={\bf 0}_{H^1} \right\}.$$ \end{proposition}
{\bf Proof :}
Let $g : {\bf h} \in L^2_a(\nu,  H^{\bf 1}_{0,0}) \to   \int_W {\bf h} \ d\nu \in H^1$. From the \textsc{Cauchy-Schwarz} inequality, $g$ is continuous, so that $F =g^{-1}(\{{\bf 0}_{H^1}\})$ is closed. On the other hand, from the very definition of $A_\nu^{0,\infty}$, it satisfies $A_\nu^{0,\infty}\subset F$. Since the set $F$ is closed, we first obtain the inclusion $cl(A_\nu^{0,\infty})\subset F$. Conversely, we invoke denseness and continuity taking averages. More accurately, assume that ${\bf h} \in L^2_a(\nu,H^{\bf 1}_{0,0})$  is such that $\int_W {\bf h} \ d\nu={\bf 0}_{H^1}$.  From Lemma 2.1. of \cite{13.}, there exists a sequence $({\bf h_n})_{n\in \mathbb{N}} \subset V^{0,\infty}_\nu $ which converges strongly to ${\bf h}$ in $L^2_a(\nu, H^{\bf 1}_{0,0})$. For all $n\in \mathbb{N}$, set ${\bf k_n}=j({\bf h_n})$, where $j$ is the map defined in the statement of Proposition~4.1, whose continuity follows from the \textsc{Cauchy-Schwarz} inequality. From Proposition~4.1, we first obtain ${\bf k_n} \in  A_\nu^{0,\infty}$. Since $({\bf h_n})_{n\in \mathbb{N}}$ converges to ${\bf h}$ and $\int_W {\bf h} \ d\nu={\bf 0}_{H^1}$, by continuity $\left(\int_W {\bf h_n} \ d\nu\right)_{n\in \mathbb{N}}$ converges to ${\bf 0}_{H^1}$ in $H^1$. Together with the triangular inequality, and with the convergence of $({\bf h_n})_{n\in \mathbb{N}}$ to ${\bf h}$, we get that $({\bf k_n})_{n\in \mathbb{N}}$ converges to ${\bf h}$. Thus, ${\bf h} \in  cl(A_\nu^{0,\infty})$. \qed

\section{Average preserving least action principle} 
\label{Section5}
\begin{lemma}\label{Lemma5.1}
Given $\nu \in {\bf M}_1(W)$, assume that $\phi : \nu \in {\bf M}_1(W) \to \phi(\nu) \in [0,+\infty]$ is $L^2_a(\nu,H^{\bf 1}_{0,0})-$ differentiable at $\nu$, and let ${\bf \xi} \in L^2_a(\nu,H^1)$ be such that  $$\delta \phi_\nu[{\bf h}] =\int_W <{\bf \xi}, {\bf h}>_{H^1}  d\nu ,$$ for all ${\bf h}\in V_\nu^{0,\infty}$. Then, we have $$\delta \phi_\nu[{\bf h}]=0, \forall {\bf h}\in L^2_a(\nu,H^1) : {\bf h_0}={\bf h_1} ={\bf 0}_{\mathbb{R}^d}, \  \nu-a.s.\  and \ \int_W {\bf h} \ d\nu ={\bf 0}_{H^1},$$ if and only if, there exist a \textit{c\`{a}d-l\`{a}g} $\mathbb{R}^d-$ valued $(\mathcal{F}_t^\nu)-$  martingale $({\bf N_t^\nu})_{t\in[0,1)}$, and a $\mathbb{R}^d-$ valued deterministic mesurable process $({\bf A_t^\nu})_{t\in [0,1)}$,  defined on the complete probability space $(W,\mathcal{B}_W^\nu,$ $\nu)$, such that $${\bf \xi} = \int_0^. {\bf A_s^\nu} ds + \int_0^. {\bf N_s^\nu} ds , \ \nu-a.s., $$ and $\int_0^1 \|{\bf A_s^\nu}\|_{\mathbb{R}^d}^2 ds < +\infty.$
\end{lemma}
{\bf Proof :}
Let ${\bf h}\in V^{0,\infty}_\nu$, and set ${\bf k}={\bf h}-\int_W {\bf h} \ d\nu$, $\nu-a.s.$. From Proposition~4.1, ${\bf k}\in A_\nu^{0,\infty}\subset V_\nu^{0,\infty}$. Moreover, we have ${\bf k_0}={\bf k_1}= {\bf 0}_{\mathbb{R}^d}$, $\nu-a.s.$, and 
$\int_W {\bf k} d\nu= {\bf 0}_{H^1}$.  Assuming that $\delta \phi_\nu[{\bf k}]= 0$, we obtain  $$ 0= \delta \phi_\nu[{\bf k}] =  \mathbb{E}_\nu\left[<{\bf \xi}, {\bf h}-\int_W {\bf h}\  d\nu>_{H^1}\right]   =  \mathbb{E}_\nu\left[<{\bf \xi}- \int_W {\bf \xi} \ d\nu, {\bf h}>_{H^1}\right]. $$ Since this holds for all ${\bf h}\in V_\nu^{0,\infty}$, from the variational characterization of martingales (see \cite{40.}, or Proposition 1.1 of \cite{13.} for a summary of the proof) and from Lemma 2.1 of \cite{13.}, we obtain the existence of a $(\mathcal{F}_t^\nu)-$martingales $({\bf N_t^\nu})_{t\in [0,1)}$ which meets the above condition, with ${\bf A_t^\nu}= \mathbb{E}_\nu[{\bf \dot{\xi}_t}]$, for all $t\in [0,1)$ outside some $\lambda_{\mathbb{R}}$-null set.  Conversely, assume the existence of such a $(\mathcal{F}_t^\nu)-$ martingale $({\bf N_t^\nu})_{t\in [0,1)}$, and of such a deterministic process $({\bf A_t^\nu})_{t\in [0,1)}$. Since $\int_0^. ({\bf N_t^\nu}-\mathbb{E}_\nu[{\bf N_0^\nu}]) dt$ is orthogonal to $L^2_a(\nu, H^{\bf 1}_{0,0})$ in $L^2_a(\nu, H^1)$ (see \cite{40.}  or Proposition 1.1. of \cite{13.}),  we obtain similarly that $$\delta \phi_\nu\left[{\bf h}-\int_W {\bf h} \ d\nu\right] =0,$$ for all ${\bf h}\in V^{0,\infty}_\nu$. Thus, from Proposition~4.1, we obtain $\delta \phi_\nu[{\bf k}]=0$ for all ${\bf k}\in A_\nu^{0,\infty}$. By Proposition~4.2, the continuity of $\delta \phi_\nu$ yields $\delta \phi_\nu[{\bf h}]= 0,$ for all ${\bf h}\in L^2_a(\nu,H^{\bf 1}_{0,0})$ such that $\int_W {\bf h} \ d\nu={\bf 0}_{H^1}.$ \qed

We now provide a recall from Definition~{5.2.} of \cite{13.}, which will be used to shorten the statement of Theorem~\ref{aplapT} below.   Given a \textsc{Borel} measurable mapping $$\mathcal{L} : (t,{\bf x},{\bf v},{\bf a}) \in[0,1]\times \mathbb{R}^d \times  \mathbb{R}^d \times ( \mathbb{R}^d\otimes  \mathbb{R}^d) \to \mathcal{L}_t({\bf x},{\bf v},{\bf a}) \in \mathbb{R}\cup\{ +\infty\},$$ and setting $$\Dom(\mathcal{L}) = \left\{ (t,{\bf x},{\bf v},{\bf a}) : \mathcal{L} < +\infty \right\},$$ $\mathcal{L}$ will be said to be a \textit{regular Lagrangian}  if it satisfies the following assumptions
\begin{enumerate}[(i)]
\item $\Dom(\mathcal{L})= [0,1]\times \mathbb{R}^d\times \mathbb{R}^d\times (\mathbb{R}^d\otimes \mathbb{R}^d)$.
\item For all $(t,{\bf x},{\bf v},{\bf a})\in \Dom(\mathcal{L})$,  $$\widetilde{\mathcal{L}}(t,{\bf x},{\bf v},{\bf a}) : ({\bf \widetilde{x}}, {\bf \widetilde{v}})\in \mathbb{R}^d\times \mathbb{R}^d  \to \mathcal{L}_t({\bf x}+{\bf \widetilde{x}},{\bf v}+ {\bf \widetilde{v}},{\bf a}) \in \mathbb{R}$$ is \textsc{Fr\'{e}chet} differentiable at ${\bf 0}_{\mathbb{R}^d\times \mathbb{R}^d}$.  
\item The mappings $(t,{\bf x},{\bf v},{\bf a})\in \Dom(\mathcal{L}) \to \bf{\partial}_{\bf x} \mathcal{L}_t({\bf x},{\bf v},{\bf a})\in \mathbb{R}^d$ and $(t,{\bf x},{\bf v},{\bf a})\in \Dom(\mathcal{L})$  $\to$ ${\bf \partial}_{\bf v} \mathcal{L}_t({\bf x},{\bf v},{\bf a})\in \mathbb{R}^d$ are \textsc{Borel} measurable.
\end{enumerate}

When $\mathcal{L}$ is a regular Lagrangian, for any $(t,{\bf x},{\bf v},{\bf a}) \in \Dom(\mathcal{L})$, we define $$D\mathcal{L}_{t,{\bf x},{\bf v},{\bf a}} :  \mathbb{R}^d 
\times \mathbb{R}^d \to  \mathbb{R}$$ by $$D\mathcal{L}_{t,{\bf x},{\bf v},{\bf a}}[{\bf \widetilde{x}},{\bf \widetilde{v}}] =  <({\bf \partial}_{\bf x} \mathcal{L}_t)({\bf x},{\bf v},{\bf a}), {\bf \widetilde{x}}>_{\mathbb{R}^d}  + <({\bf \partial}_{\bf v} \mathcal{L}_t)({\bf x},{\bf v},{\bf a}), {\bf \widetilde{v}}>_{\mathbb{R}^d},$$ the linear operator such that $$D\mathcal{L}_{t,{\bf x},{\bf v},{\bf a}}[{\bf \widetilde{x}}, {\bf \widetilde{v}}] = \frac{d}{d\epsilon}\mathcal{L}_t({\bf x}+ \epsilon {\bf \widetilde{x}} ,{\bf v} + \epsilon {\bf \widetilde{v}}, {\bf a} )\big|_{\epsilon =0},$$ for any $({\bf \widetilde{x} }, {\bf \widetilde{v}})\in \mathbb{R}^d\times \mathbb{R}^d$, from which the notation ${\bf \partial}_{\bf x} \mathcal{L}_t$ (respectively ${\bf \partial}_{\bf v} \mathcal{L}_t$) is clear.

\begin{theorem} \label{aplapT}
Let $\mathcal{L}$ be a regular Lagrangian whose associated action functional on ${\bf M}_1(W)$ is defined by

$$ \phi(\nu) = \begin{cases}
        \mathbb{E}_\nu\left[\int_0^1 \mathcal{L}_t({\bf W_t},{\bf v_t^\nu},{\bf \alpha_t^\nu})dt \right]  & \text{if } \   \mathbb{E}_\nu\left[\int_0^1 |\mathcal{L}_t({\bf W_t}, {\bf v_t^\nu},{\bf \alpha_t^\nu})|dt \right] < +\infty \\
       + \infty & \text{otherwise}
       \end{cases}, $$ for all $\nu \in \bf{M_1^{\mathbb{S}}}(W)$, and by   $\phi(\nu) =+\infty$ if $\nu \in \bf{M_1}(W)\backslash \bf{M_1^{\mathbb{S}}}(W).$
 Further assume  the existence of a strictly positive continuous function $f :\mathbb{R}^d \to \mathbb{R}^+$ and of $p_1,p_2\geq 2$ such that
$$ \limsup_{|\epsilon| \downarrow 0} \sup_{(t, {\bf x}, {\bf v},{\bf a}, {\bf \widetilde{x}}, {\bf \widetilde{v}}) \in \Dom(\mathcal{L}) \times \mathbb{R}^d\times \mathbb{R}^d} F_\epsilon(t,{\bf x},{\bf v},{\bf a},{\bf \widetilde{x}},{\bf \widetilde{v}})=0, $$ where $$F_\epsilon(t, {\bf x}, {\bf v},{\bf a}, {\bf \widetilde{x}}, {\bf \widetilde{v}})= \frac{\left| \mathcal{L}_t(x+\epsilon {\bf \widetilde{x}}, {\bf v}+ \epsilon{\bf \widetilde{v}},{\bf a}) -  \mathcal{L}_t({\bf x},{\bf v},{\bf a}) - \epsilon D\mathcal{L}_{t,{\bf x},{\bf v},{\bf a}}[{\bf \widetilde{x}},{\bf \widetilde{v}}]  \right|}{\epsilon f({\bf \widetilde{x}})\left(1+ \|{\bf \widetilde{v}}\|_{\mathbb{R}^d}^2 + G(t,{\bf x},{\bf v},{\bf a})  \right)},$$ for all $(\epsilon,t,{\bf x},{\bf v},{\bf a},{\bf \widetilde{x}},{\bf \widetilde{v}})\in \mathbb{R}\times \Dom(\mathcal{L}) \times \mathbb{R}^d\times \mathbb{R}^d$, and where $$G(t,{\bf x},{\bf v},{\bf a})= |\mathcal{L}_t({\bf x},{\bf v},{\bf a})| + \|{\bf \partial}_{\bf x} \mathcal{L}_t({\bf x},{\bf v},{\bf a})\|_{\mathbb{R}^d}^{p_1}  + \|{\bf \partial}_{\bf v} \mathcal{L}_t({\bf x},{\bf v},{\bf a})\|_{\mathbb{R}^d }^{p_2},$$ for all $(t,{\bf x},{\bf v},{\bf a})\in \Dom(\mathcal{L})$. Then, for any $\nu \in {\bf M}_1^{\mathbb{S}}(W)$ which meets the integrability condition  $$ \phi(\nu) + \mathbb{E}_\nu\left[\int_0^1 \left(\|{\bf \partial}_{\bf x} \mathcal{L}_s({\bf W_s}, {\bf v_s^\nu}, {\bf \alpha_s^\nu})\|^{p_1}_{\mathbb{R}^d} + \|{\bf \partial}_{\bf v} \mathcal{L}_s({\bf W_s},{\bf v_s^\nu},{\bf \alpha_s^\nu})\|^{p_2}_{\mathbb{R}^d}\right) ds \right] < +\infty,$$ we have that $\phi$ is $L^2_a(\nu,H^{\bf 1}_{0,0})-$differentiable at $\nu$. Moreover, in this case, the following assertions are equivalent

 \begin{enumerate}[(i)]
 \item  We have  $\delta \phi_\nu[{\bf h}] =0$,  $\forall {\bf h}\in L^2_a(\nu, H^1)$ such that  ${\bf h_0}={\bf h_1}={\bf 0}_{\mathbb{R}^d}, \ \nu-a.s.$,  and $\int_W {\bf h} \ d\nu = {\bf 0}_{H^1}.$ 
 \item   $\nu$ satisfies the following \textsc{Euler-Lagrange} condition : there exists a \textit{c\`{a}d-l\`{a}g} $\mathbb{R}^d-$ valued $(\mathcal{F}_t^\nu)-$ martingale $({\bf N_t^\nu})
 _{t\in[0,1)}$, and a deterministic measurable process $({\bf A_t^\nu})_{t\in[0,1)}$, defined on the complete probability space $(W,\mathcal{B}_W^\nu,\nu)$, such that \begin{equation} \label{ELA} {\bf \partial}_{\bf v} 
 \mathcal{L}_t({\bf \omega(t)},{\bf v_t^\nu}(\omega),{\bf \alpha_t^\nu}(\omega)) - \int_0^t {\bf \partial}_{\bf x} \mathcal{L}_s({\bf \omega(s)},{\bf v_s^\nu}(\omega),{\bf \alpha_s^\nu}(\omega)) ds  = {\bf A_t^\nu}+ {\bf N^\nu_t}(\omega), \end
 {equation} holds for all $(\omega, t) \in W \times [0,1)$, outside a $\nu\otimes\lambda_{\mathbb{R}}$-null set.  Moreover, we have $$\int_0^1\|{\bf A_s^\nu}\|^2_{\mathbb{R}^d}ds + 
 \mathbb{E}_\nu\left[\int_0^1 \|{\bf N_s^\nu}\|_{\mathbb{R}^d}^2 ds \right] < +\infty.$$
\end{enumerate}
  \end{theorem}
{\bf Proof :}
For $t\in[0,1)$, define \begin{equation} \label{eqrevfrev1} {\bf \dot{\xi}_t}  = {\bf \partial}_{\bf v}\mathcal{L}_t({\bf W_t},{\bf v_t^\nu},{\bf \alpha_t^\nu}) -\int_0^t {\bf \partial}_{\bf x} \mathcal{L}_s({\bf W_s},{\bf v_s^\nu}, {\bf \alpha_s^\nu}) ds, \end{equation} and notice that, from~(\ref{eqrevfrev1}), the integrability condition ensures that, for ${\bf \omega}\in W$ outside a specific $\nu-$negligible set, we have $\int_0^1 |{\bf \dot{\xi}(\omega)}|^2 ds <+\infty$, so that we also have $\int_0^1 |{\bf \dot{\xi}(\omega)}| ds <+\infty$. Define ${\bf \xi}=  \int_0^. {\bf \dot{\xi}_t} dt, \ \nu-a.s..$  Under those conditions, the $L^2_a(\nu,H^{\bf 1}_{0,0})-$ differentiability of $\phi$ follows from Theorem 5.1 of \cite{13.}.  Moreover, the proof of the latter also yields  $$\delta \phi_\nu[{\bf h}]= \int_W <{\bf \xi}, {\bf h}>_{H^1} d\nu,$$ for all ${\bf h}\in L^2_a(\nu,H^{\bf 1}_{0,0})$. Whence, by applying Lemma~\ref{Lemma5.1}, together with classical methods (for instance see Lemme VIII.1 of \cite{41.}), the result follows.   \qed

 \section{Forward-Backward systems of \textsc{Mckean-Vlasov} stochastic equations with classical actions}
 \label{Section6}
 
In this section we take a classical action, with the usual convention of \textbf{E}uclidean \textbf{Q}uantum \textbf{M}echanics (\textbf{EQM}) on Lagrangians (for more details on \textbf{EQM}, see \cite{42.}, \cite{43.} and \cite{19.}); for the sake of clarity, the function $\mathcal{V} : \mathbb{R}^d\to \mathbb{R}$ will be assumed to be smooth. This section involves systems of stochastic differential equations based on \textsc{It\^{o}}'s stochastic integrals (see \cite{45.} and \cite{44.}), and the so-called driving process is a \textit{Brownian motion}, whose law on $W$, as a random path or function, is the so-called classical \textsc{Wiener} measure (see \cite{46.}).   Subsequently, we use the usual notation $$\mathbb{E}_{\mathcal{P}}\left[{\bf X}\right]=   \int_{\mathbb{R}^d} {\bf x}  \ X_{\star}\mathcal{P}({\bf dx}) = \sum_{i=1}^d \mathbb{E}_{\mathcal{P}}\left[<{\bf X},{\bf e_i}>_{\mathbb{R}^d}\right] {\bf e_i},$$ $({\bf e_i})_{i\in \{1,..., d\}}$ denoting the canonical basis of $\mathbb{R}^d$, where ${\bf X} : \Omega \to \mathbb{R}^d$ is a $\mathcal{P}-$ integrable measurable function on a complete probability space $(\Omega, \mathcal{A}, \mathcal{P})$.

\begin{proposition} \label{PropMVFB}  \label{FBPROPMV} Let $d\in \mathbb{N}$, $d\geq 1$. Given two \textsc{Borel} probability measures  $\nu_0, \nu_1 \in {\bf M}_1(\mathbb{R}^d)$, let $$\mathcal{L}^V_t({\bf x},{\bf v},{\bf a})=  \frac{\|{\bf v}\|_{\mathbb{R}^d}^2}{2} +\mathcal{V}({\bf x}) ,$$  for all $(t,{\bf x},{\bf v},{\bf a})\in [0,1]\times \mathbb{R}^d\times \mathbb{R}^d\times(\mathbb{R}^d\otimes \mathbb{R}^d)$, where $\mathcal{V}\in \mathcal{C}^\infty(\mathbb{R}^d; \mathbb{R})$, and let $({\bf \sigma_t})_{t\in[0,1]}$ be a $\mathbb{R}^d\otimes \mathbb{R}^d-$ valued predictable process (see \cite{29.}) on the canonical space $(W,\mathcal{B}_W)$; ${\bf \sigma} : (t,\omega) \in [0,1]\times W\to \sigma_t(\omega) \in \mathbb{R}^d\otimes \mathbb{R}^d$. Then, the following assertions are equivalent :

\begin{enumerate}[(i)]

\item There exist a \textit{complete stochastic basis} $(\Omega, \mathcal{A}, (\mathcal{A}_t),$ $\mathcal{P})$, an $(\mathcal{A}_t)-$Brownian motion $({\bf B_t})_{t\in[0,1]}$, a \textit{c\`{a}d-l\`{a}g} $(\mathcal{A}_t)-$martingale $({\bf Z_t})_{t\in[0,1)}$ on this space, and a pair of measurable $(\mathcal{A}_t)-$ adapted processes $({\bf X_t})$, $({\bf Y_t})$, where $({\bf X_t})_{t\in[0,1]}$ is a continuous process, and where $({\bf Y_t})_{t\in[0,1)}$ is a \textit{c\`{a}d-l\`{a}g} process,  which solve the following system :
 \begin{equation*} \begin{cases}{\bf X_t} = {\bf X_0} + \int_0^t {\bf \sigma_s}({\bf X}) d{\bf B_s} + \int_0^t {\bf Y_s} ds  \\  {\bf Y_t}=   {\bf Z_t} + \int_{\mathbb{R}^d} {\bf y} \ {{\bf Y_t}_\star \mathcal{P}}({\bf d y}) + \int_0^t \left( \int_{\mathbb{R}^d} \left( {\bf \nabla \mathcal{V}}({\bf X_s})- {\bf \nabla \mathcal{V}}({\bf x}) \right) \ {{\bf X_s}}_\star \mathcal{P}({\bf dx}) \right) ds, 
 \\ {{\bf X_0}}_\star \mathcal{P} = \nu_0, \   {{\bf X_1}}_\star \mathcal{P} = \nu_1 
 \end{cases}  \end{equation*}  together with the integrability conditions 
  \begin{equation} \mathbb{E}_{\mathcal{P}}\left[\int_0^1 \|{\bf Y_s}\|^2_{\mathbb{R}^d} ds\right]+ \mathbb{E}_{\mathcal{P}}\left[\int_0^1 \|{\bf \nabla 
\mathcal{V}}({\bf X_s})\|^2_{\mathbb{R}^d} ds \right] < +\infty, \label{intloc} \end{equation} and
$$ \sum_{i=1}^d \sum_{j=1}^d \int_0^1 |({\bf \sigma_s}.{\bf \sigma_s}^\dagger)^{i,j}({\bf X})| ds < +\infty ,  \ \mathcal{P}-
a.s.$$ \item There exists $\nu \in {\bf M}_1^{\mathbb{S}}(W)$, with  ${{\bf W_0}}_\star \nu = \nu_0$,  and ${{\bf W_1}}_\star \nu = \nu_1,$ which is such that~(\ref{ELA}) holds with $\mathcal{L}^V$, for some c\`{a}d-l\`{a}g $(\mathcal{F}_t^\nu)-
$martingale $({\bf N_t^\nu})_{t\in[0,1)}$, and for some c\`{a}d-l\`{a}g deterministic process $({\bf A^\nu_t})_{t\in[0,1)}$, both defined on the complete probability space $(W, \mathcal{B}_W^\nu, \nu)$. Moreover, $\nu$ further satisfies  $$\mathbb{E}_\nu\left[\int_0^1 \left(\|{\bf v^\nu_s}\|^2_{\mathbb{R}^d} +  \|{\bf \nabla \mathcal{V}}({\bf W_s})\|^2_{\mathbb{R}^d}\right) ds \right] < +\infty,$$  and $$\int_0^. {\bf \alpha_t^\nu} dt = \int_0^.({\bf \sigma_t}. {\bf \sigma_t^\dagger}) dt, \ \nu-a.s..$$ 
\end{enumerate} 
Finally, in this case we can take $\nu = X_{\star}\mathcal{P}$. 
\end{proposition}
{\bf Proof :}
 First assume the existence of a pair of measurable processes $({\bf X},{\bf Y})$, where $({\bf X_t})_{t\in[0,1]}$ is continuous and $({\bf Y_t})_{t\in[0,1)}$ is \textit{c\`{a}d-l\`{a}g},  defined on a \textit{complete stochastic basis} $(\Omega, \mathcal{A},$ $(\mathcal{A}_t), \mathcal{P})$, and $(\mathcal{A}_t)-$ adapted, as in the statement $(i)$. Let $({\bf Z_t})_{t\in[0,1)}$ be a $\mathbb{R}^d-$valued \textit{c\`{a}d-l\`{a}g}  $(\mathcal{A}_t)-$ martingale which satisfies the second stochastic differential equation of the system in $(i)$. From standard results on transformations of laws of semi-martingales (for instance see \cite{13.}), we first 
obtain $\nu={\bf X}_\star \mathcal{P}\in {\bf M}_1^{\mathbb{S}}(W)$. We also obtain the conditions on \textit{local characteristics} \begin{equation} \int_0^. <{\bf v_t^\nu} \circ {\bf X},{\bf e_i}>_{\mathbb{R}^d} dt = \int_0^.\mathbb{E}_{\mathcal{P}}\left[Y^i_t | \mathcal{G}_t^{\bf X}\right] dt, \ \mathcal{P}-a.s., \label{prop51f1} \end{equation} and $$\int_0^. (\alpha_s^\nu)^{i,j}\circ {\bf X} ds = \int_0^. ({\bf \sigma_s}. {\bf \sigma_s}^\dagger)^{i,j} ds, \  \mathcal{P}-a.s., $$ for all $i,j\in\{1,...,d\}$, where $(\mathcal{G}_t^{\bf X})$ denotes the $\mathcal
{P}-$usual augmentation of the filtration generated by ${\bf X}$. For $t\in[0,1)$, let ${\bf Z}^t_s(\omega)={\bf Z}_{s\wedge t}(\omega)$, $\forall \omega\in \Omega$, and $s\in[0,1)$; in 
particular, $({\bf Z^t_s})_{s\in[0,1)}$ is an $(\mathcal{A}_s)-$ martingale, so that $(\|{\bf Z_s^t}\|^2_{\mathbb{R}^d})_{s\in[0,1)}$ is an $(\mathcal{A}_s)-$ submartingale. Therefore, the function $s\in [0,1) \to \mathbb{E}_{\mathcal{P}}[\|{\bf Z_s^t}\|^2_{\mathbb{R}^d}] \in \mathbb{R}_+$ is increasing.  Together with~(\ref{intloc}) and with \textsc{Jensen}'s inequality, it implies  that \begin{equation} \label{ZmajorMV} \sup_{s\in [0,1)} \mathbb{E}_{\mathcal{P}}\left[\|Z_s^t\|^2_{\mathbb{R}^d}\right] <+\infty.\end{equation} Hence, for $i\in \{1,...,d\}$, the \textsc{La Vall\'{e}e-Poussin} criterion (see \cite{29.}) ensures that $(Z_s^t)_{s\in [0,1)}$ is a uniformly integrable martingale, and in 
particular that it is a martingale of class $D$ (see \cite{30.}). Therefore, according to \cite{BAIN} p.19, thanks to a classical theorem (see 47, p.119 of \cite{30.}), for $i\in \{1,...,d\}$, we are allowed to conclude that the optional projection $(\widehat{Z}_s^{t,i})_{s\in [0,1)}$ of $(<{\bf Z^t_s}, {\bf e_i}>_{\mathbb{R}^d})_{s\in[0,1)}$ on $(\mathcal{G}_t^{\bf X})$ is \textit{c\`{a}d-l\`{a}g}.  Then, it is enough to notice that, from the definitions and using the right-continuity, we have $\widehat{Z}^{t,i}_{u\wedge t} = \widehat{Z}^{\widetilde{t},i}_{u\wedge t}$, $\forall u\in [0,1)$ , $\mathcal{P}-a.s.$, whenever $0\leq t \leq \widetilde{t}<1$; hence, we may define $\widehat{Z}^{i}_u =\lim_{t\to 1} \widehat{Z}^{t,i}_u$, $\forall u\in [0,1)$, $\mathcal{P}-a.s.$. From its construction, $(\widehat{Z}^{i}_u)_{u\in[0,1)}$ is a \textit{c\`{a}d-l\`{a}g} modification of $(E_{\mathcal{P}}[<{\bf Z_u},{\bf e_i}>_{\mathbb{R}^d}| \mathcal{G}_u^{\bf X}])_{u\in[0,1)}$ ; henceforth and until the end of the proof, we denote by $(\mathbb{E}_{\mathcal{P}}\left[Z_t^i| \mathcal{G}_t^{\bf X}\right])$ for $(\widehat{Z}^{i}_t)$. Then, a \textit{c\`{a}d-l\`{a}g} modification of $(\mathbb{E}_{\mathcal{P}}\left[Y_t^i | \mathcal{G}_t^{\bf X}\right] )$  is obtained by setting 
\begin{equation} \mathbb{E}_{\mathcal{P}}\left[Y_t^i  | \mathcal{G}_t^{\bf X}\right]=   \mathbb{E}_{\mathcal{P}}[Y^i_t] + \mathbb{E}_{\mathcal{P}}\left[Z_t^i| \mathcal{G}_t^{\bf X}\right] +\int_0^t \left(\frac{\partial\mathcal{V}}{\partial x^i}({\bf X_s})- \mathbb{E}_{{{\bf X_s}}_
\star \mathcal{P}}\left[\frac{\partial \mathcal{V}}{\partial x^i}\right]\right) ds, \label{prop51f2}\end{equation} for all $t\in[0,1)$, and $i\in\{1,...,d\}$, $\mathcal{P}-a.s.$; notice that the hypothesis on $({\bf Z_t})$ and $({\bf Y_t})$, notably the equation satisfied by $({\bf Y_t})$,  yield that  the function $t\in [0,1) \to \mathbb{E}_{\mathcal{P}}[Y^i_t]\in \mathbb{R}$ is  \textit{c\`{a}d-l\`{a}g}. Moreover, for all $u\in [0,1)$, from the process whose components are $(\mathbb{E}_{\mathcal{P}}\left[Y_t^i | \mathcal{G}_t^{\bf X}\right]$ $)_{t\in[0,u]}$, $i\in\{1,...,d\}$,  we can obtain an $\mathcal{A}/ \mathcal{B}_{D[0,u]}-$measurable 
function ${\bf F}^u : \Omega \to D([0,u])$, where $D_u=D([0,u])$, the space of $\mathbb{R}^d-$ valued \textit{c\`{a}d-l\`{a}g} functions defined on the closed interval $[0,u]$, is endowed with the \textsc{Skorokhod} 
topology (see \cite{47.}), which turns it into a Polish space.   Whence, since $({\bf X}_t)_{t\in[0,u]}$ can be seen as a 
$W_u=\mathcal{C}([0,u], \mathbb{R}^d)$-valued measurable map ${\bf X}^u$, since the filtration $(\mathcal{A}_t)$ satisfies the usual conditions, and since $({\bf F_t^u})$ is $(\mathcal{G}_t^
{\bf X})-$adapted, where ${\bf F_t^u}= {\bf W_t}\circ {\bf F}^u$ for all $t\in [0,u]$ (see  Proposition 1.3 of  \cite{13.}), we 
can find a $\mathcal{B}_{W_u}^{\nu_u}/ \mathcal{B}_{D[0,u]}-$measurable function ${\bf v}^u : W_u\to D_u$ such that ${\bf F}^u= {\bf v}^u\circ {\bf X}^u$, $\mathcal{P}-a.s.$, and $({\bf v_t^u})_{t\in [0,u]}$ is $(\mathcal{F}_t^{\nu_u})_{t\in [0,u]}-$ adapted, $(\mathcal{F}_t^{\nu_u})_{t\in [0,u]}$ denoting the $\nu_u$ usual augmentation of the filtration generated by the evaluation process on $(W_u, \mathcal{B}_{W_u}^{\nu_u})$; $\nu_u$ denotes the law $j^u_\star \nu$, where $j^u : W\to W_u$ denotes the continuous function which is defined by ${\bf j^u(\omega)(t)} ={\bf \omega(t)}$, $\forall t\in [0,u]$, $\forall \omega\in W$, while we have ${\bf X^u}= j^u\circ {\bf X}$, $\mathcal{P}-$ a.s.. From the construction, we obtain ${\bf v}_t^u\circ j^u= {\bf v}_{t}^{\widetilde{u}}\circ j^{\widetilde{u}}$, $\nu-a.s.$, for all $0\leq t\leq u \leq \widetilde{u}$, and both are \textit{c\`{a}d-l\`{a}g} on $[0,u]$. Therefore, we can define ${\bf v}_t(\omega)= \lim_{u\to 1} {\bf v}_t^u\circ j^u(\omega)$, for all $t\in[0,1)$, and for any $\omega\in W$ outside a $\nu-$negligible set. By construction, $({\bf v}_t)_{t\in[0,1)}$ is a $(\mathcal{F}_t^\nu)-$ adapted \textit{c\`{a}d-l\`{a}g} process. Moreover, the definitions of ${\bf F}^u$, $\forall u\in [0,1)$, and of $({\bf v}_t)_{t\in[0,1)}$, yield  $$<{\bf v_t}\circ {\bf X}, {\bf e_i}>_{\mathbb{R}^d} = \mathbb{E}_{\mathcal{P}}\left[Y_t^i | \mathcal{G}_t^{\bf X}\right] ,$$ for all $t\in[0,1)$, $i\in\{1,...,d\}$, $\mathcal{P}-a.s.$. From the above condition~(\ref{prop51f1}) on \textit{local characteristics}, we see that we can obtain a version of $({\bf v_t^\nu})$ from $({\bf v_t})$ (indeed $\int_0^. {\bf v_s} ds = \int_0^. {\bf v_s^\nu} ds$, $\nu-a.s.$). Thus, until the end of the proof we take ${\bf v_t^\nu}
({\bf \omega})= {\bf v_t}({\bf \omega}),$ for all ${\bf \omega} \in W$, and $t\in [0,1)$. From this,  $({\bf v_t^\nu})$ can be taken to be \textit{c\`{a}d-l\`{a}g}. Moreover, using the tower property, we obtain $$\mathbb{E}_{\mathcal{P}}[Y_t^i] = \mathbb{E}_{\mathcal{P}}\left[\mathbb{E}_{\mathcal{P}}\left[Y_t^i |\mathcal{G}_t^{\bf X}\right] \right] = \mathbb{E}_{\nu}\left[<{\bf v_t^\nu},{\bf e_i}>_{\mathbb{R}^d}\right],$$ for all $t\in[0,1)$, $i\in \{1,...,d\}$. Since $({\bf \widehat{Z}_t})_{t\in [0,1)}$ is \textit{c\`{a}d-l\`{a}g}, similarly as above,   we obtain a \textit{c\`{a}d-l\`{a}g} $(\mathcal{F}_t^\nu)-
$ adapted process $({\bf \widetilde{Z}_t})_{t\in [0,1)}$ on $(W,\mathcal{B}_W^\nu,\nu)$ such that $\mathbb{E}_{\mathcal{P}}\left[Z_t^i| \mathcal{G}_t^{\bf X}\right] =\widetilde{Z}_t^i\circ {\bf X},$ for all $t\in[0,1)$, $i\in\{1,...,d\}$, $\mathcal{P}-a.s..$ Setting  $$c= \mathbb{E}_\nu\left[(<{\bf \widetilde{Z}_t},{\bf e_i}>_{\mathbb{R}^d}-
<{\bf \widetilde{Z}_s}, {\bf e_i}>_{\mathbb{R}^d})1_A\right],$$ we obtain $$c =\mathbb{E}_{\mathcal{P}}\left[\left(\mathbb{E}_{\mathcal{P}}\left[Z_t^i | \mathcal{G}_t^{\bf X}\right]- \mathbb{E}_{\mathcal{P}}\left[Z_s^i | \mathcal{G}_s^{\bf X}\right] \right) {\bf 1}_{{\bf X}^{-1}(A)} \right],$$ for any $i\in \{1,...,d\}$, $0\leq s\leq t<1$ and $A\in \mathcal{F}_s^\nu$. 
Since $({\bf X_t})$ is continuous, ${\bf X}^{-1}(A)\in \mathcal{G}_s^{\bf X}$; recall that ${\bf X} : \Omega \to W$. Therefore, we get $$\mathbb{E}_\nu\left[(<{\bf \widetilde{Z}_t}, {\bf e_i}>_{\mathbb{R}^d}-<{\bf \widetilde{Z}_s}, {\bf e_i}>_{\mathbb{R}^d})1_A\right] =\mathbb{E}_{\mathcal{P}}\left[ (Z^i_t - Z^i_s)  {\bf 1}_
{{\bf X}^{-1}(A)} \right],$$ for all $i\in\{1,...,d\}$. Together with  $\mathbb{E}_{\mathcal{P}}\left[Z^i_t- Z^i_s | \mathcal{A}_s\right] = 0,$ and $\mathcal{G}_s^{\bf X} \subset \mathcal{A}_s$, for all $0\leq s\leq t <1$,  the tower property yields $$\mathbb{E}_\nu\left[(<{\bf \widetilde{Z}_t},{\bf e_i}>_{\mathbb{R}^d}-<{\bf \widetilde{Z}_s},{\bf e_i}>_{\mathbb{R}^d})1_A\right]  = 0,$$ for all $A\in \mathcal{F}_s^\nu$, $i\in\{1,...,d\}$. Since $({\bf \widetilde{Z}_t})_{t\in[0,1)}$ is also $(\mathcal{F}_t^\nu)-$adapted, we obtain that it is a $(\mathcal{F}_t^\nu)-$ martingale. Whence, using both the above condition on \textit{local characteristics}, and~on the optional projections, for $i\in \{1,...,d\}$, from~(\ref{prop51f2}) we obtain
$$<{\bf v_t^\nu}\circ {\bf X},{\bf e_i}>_{\mathbb{R}^d} - \mathbb{E}_\nu[<{\bf v_t^\nu},{\bf e_i}>_{\mathbb{R}^d}] = \widetilde{Z}^i_t\circ X + \int_0^t \left(\frac{\partial\mathcal{V}}{\partial x^i}(X_s)- \mathbb{E}_{{X_s}_\star \mathcal{P}}\left[\frac{\partial \mathcal{V}}{\partial x^i}\right]\right) ds,$$ $\forall t\in[0,1)$, $\mathcal{P}-a.s.$ for any $i\in\{1,...,d\}$, so that  $${\bf v_t^\nu} -\int_0^t {\bf \nabla \mathcal{V}}({\bf X_s}) ds = {\bf \widetilde{Z}_t} + {\bf \widetilde{A}_t}, \ \forall t\in [0,1), \ \nu-a.s.,$$ where $$\widetilde{A}^i_t= \mathbb{E}_\nu[<{\bf v_t^\nu},{\bf e_i}>_{\mathbb{R}^d}] - \int_0^t \mathbb{E}_{{{\bf W_s}}_\star\nu}\left
[\frac{\partial \mathcal{V}}{\partial x^i} \right] ds,$$ $\forall t\in [0,1)$. From the definitions together with~(\ref{intloc}) and~(\ref{ZmajorMV}), for $i\in \{1,...,d\}$, we obtain $$\sup_{t\in [0,u]}\mathbb{E}_{\nu}\left[(<{\bf v_t^\nu},{\bf e_i}>_{\mathbb{R}^d})^2 \right]<+\infty,$$ for all $u\in [0,1)$. Therefore, still from the \textsc{La Vall\'{e}e-Poussin} criterion, we deduce that the family of random variables $\{ <{\bf v_t^\nu},{\bf e_i}>_{\mathbb{R}^d} : t\in [0,u] \}$ is uniformly integrable, for all $u\in [0,1)$. Therefore, we are allowed to apply the \textsc{Lebesgue} convergence theorem of \cite{29.}, to conclude that the function $t\in [0,1) \to  \mathbb{E}_\nu[<{\bf v_t^\nu},{\bf e_i}>_{\mathbb{R}^d}] \in \mathbb{R}$ is \textit{c\`{a}d-l\`{a}g}, since the stochastic process $({\bf v_t^\nu})$ is \textit{c\`{a}d-l\`{a}g}. As a consequence,  $(\widetilde{A}_t^i)_{t\in[0,1)}$ is also \textit{c\`{a}d-l\`{a}g}, $\forall i\in\{1,...,d\}$. Conversely, the result follows easily by taking a good extension of the canonical space if necessary, and by applying the martingale representation 
theorem of \cite{45.}. \qed

\begin{remark}
From the above proof, the martingale $({\bf Z_t})_{t\in[0,1)}$ in $(i)$ of the statement of Proposition~\ref{PropMVFB} is a \textit{locally square integrable} $(\mathcal{A}_t)-$ martingale (see \cite{PROTTER}). (End of Remark.)
\end{remark}

   \begin{example}
Let $f\in C([0,1],\mathbb{R}^d)$ and $b\in C^{1,2}([0,1]\times \mathbb{R}^d, \mathbb{R}^d)$ be differentiable (respectively twice differentiable) in the first (respectively second) variable with continuous partial derivatives, which further satisfiy 
$$\frac{\partial {\bf b}}{\partial t} + ({\bf b}. \nabla)  {\bf b} + \frac{\Delta {\bf b} }{2} = {\bf f},$$
and the \textsc{Novikov} condition (see \cite{45.}, \cite{PROTTER}) $$\int_W \exp\left(\frac{1}{2}\int_0^1\| {\bf b}(t,{\bf \omega}(t))\|^2_{\mathbb{R}^d} dt\right) \mu_W(d{\bf \omega}) <+\infty,$$ $\mu_W$ denoting the classical \textsc{Wiener} measure (see \cite{46.}) on the measurable space $(W,$ $\mathcal{B}_W)$. As a consequence, $$\frac{d\nu}{d\mu_W}= \exp\left(\int_0^1 <{\bf b}(t,{\bf W_t}), d{\bf W_t}> - \frac{1}{2}\int_0^1 \|{\bf b}(t,{\bf W_t})\|_{\mathbb{R}^d}^2 dt \right), \ \mu_W-a.s.$$ defines the \textsc{Radon-Nikodym} derivative of a \textsc{Borel} probability measure $\nu$, which is absolutely continuous with respect to $\mu_W$. Further assume that \begin{equation} \int_{[0,1]\times W}\|{\bf b}(t, {\bf \omega}(t))\|_{\mathbb{R}^d}^2 \lambda\otimes \nu(dt, {\bf d\omega})<+\infty. \label{girscondexa1}\end{equation} Together's with \textsc{L\'{e}vy}'s criterion (see \cite{45.}, \cite{PROTTER}), the \textsc{Girsanov} theorem (see \cite{GIRSANOV}, \cite{PROTTER}) yields $\nu \in {\bf M}_1^{\mathbb{S}}(W)$, with ${\bf \alpha_t^\nu}({\bf \omega}) = I_{\mathbb{R}^d}$, and ${\bf v_t^\nu}({\bf \omega})= b(t,{\bf \omega(t)})$, $\nu\otimes \lambda|_{[0,1]}-a.e..$  Since~(\ref{girscondexa1})  holds, by applying \textsc{It\^{o}}'s formula (see \cite{45.}), we obtain that $\nu$ satisfies~(\ref{ELA}), with $$\mathcal{L}_t({\bf x},{\bf v},{\bf a})= \frac{\|{\bf v}\|_{\mathbb{R}^d}^2}{2},$$ for any $(t,{\bf x},{\bf v}, \bf{a})$ $\in$ $[0,1]\times \mathbb{R}^d\times \mathbb{R}^d\times (\mathbb{R}^d\otimes \mathbb{R}^d)$.
\end{example}
   
  \begin{example} \label{ex2}(The Brownian motion with a \textsc{L\'{e}vy} drift.) 
  For $d\in \mathbb{N}$, $d\geq 1$, let $({\bf B_t})_{t\in[0,1]}$ be an $\mathbb{R}^d-$ valued $(\mathcal{A}_t)-$ Brownian motion defined on a complete stochastic basis $(\Omega, \mathcal{A},$ $(\mathcal{A}_t)_{t\in[0,1]}, \mathcal{P})$, and let $({\bf Y_t})_{t\in [0,1)}$ be a $\mathbb{R}^d-$ valued $(\mathcal{A}_t)-$\textsc{L\'{e}vy} 
  process (see \cite{45.}, \cite{KunitaJump}, \cite{PROTTER}, \cite{SATO}), which is defined on the same stochastic basis, with characteristic \textsc{L\'{e}vy} triplet $({\bf A},{\bf b}, m)$, where ${\bf A} = {\bf \sigma}.{\bf \sigma}^t$ for some $d\times \widetilde{d}$ matrix ${\bf \sigma}$ whose rank is  $\widetilde{d}\leq d$ (see \cite{45.}), where ${\bf b}\in 
  \mathbb{R}^d$, and where $m$ denotes the so-called \textsc{L\'{e}vy} measure of $(Y_t)$, which is a $\sigma-$finite measure on $\mathbb{R}^d\setminus \{{\bf 0_{\mathbb{R}^d}}\}$, such that $$\int_
  {\mathbb{R}^d\setminus\{{\bf 0_{\mathbb{R}^d}}\}} \min\left(\|{\bf x}\|^2_{\mathbb{R}^d},1\right)  m({\bf dx}) <+\infty.$$ Recall that, according to the \textsc{L\'{e}vy-It\^{o}} theorem (see \cite{45.}, \cite{KunitaJump}, \cite{SATO}, \cite{StroockPTAV}), on a possibly enlarged probability space, $({\bf Y_t})$ can be represented 
  under the form  $$ <{\bf Y_t},{\bf e_i}>_{\mathbb{R}^d} = <{\bf Y_0},{\bf e_i}>_{\mathbb{R}^d} + \sum_{j=1}^{\widetilde{d}} {\bf \sigma}_{ij} \widetilde{B}^j_t + <{\bf b},{\bf e_i}>_{\mathbb{R}^d} t + I_t^i + J_t^i,$$
   where $$I_t^i =   \int_{0}^{t}\int_{\mathbb{R}^d\setminus {B}({\bf 0_{\mathbb{R}^d}}, 1)} <{\bf x},{\bf e_i}>_{\mathbb{R}^d} N(ds,{\bf dx}),$$   ${B}({\bf 0_{\mathbb{R}^d}}, 1)$ denoting the unit ball of $\mathbb{R}^d$ endowed with its Euclidean norm,  where $$J_t^i= \int_{0}^{t}\int_{{B}({\bf 0_{\mathbb{R}^d}}, 1) \setminus\{{\bf 0_{\mathbb{R}^d}}\}} <{\bf x},{\bf e_i}>_{\mathbb{R}^d} \widetilde{N}(ds,{\bf dx}),$$ for all $t\in [0,1)$, almost surely, $\forall i\in \{1,...,d\}$, and where $( \widetilde{B}_t)$ 
  is a $\mathbb{R}^{\widetilde{d}}$- valued Brownian motion, while  $N$ denotes a \textsc{Poisson} random measure with intensity measure $\lambda\otimes m$ on $[0,1)\times (\mathbb{R}^d\setminus \{{\bf 0_{\mathbb{R}^d}}\})$, which is independent of $(\widetilde{B}_t)$, and whose  compensated random measure $\widetilde{N}$ is given by  $\widetilde{N} = N- \lambda\otimes m$;  the term with integral of $\widetilde{N}$ can be constructed by a $L^2$ convergence from steps functions (see \cite{KunitaJump}). Moreover,  for $t\in (0,1)$,  the probability measure $({\bf Y_{t}}-{\bf Y}_0) _\star\mathcal{P}\in {\bf M}_1(\mathbb{R}^d)$ is 
  infinitely divisible, with characteristic function provided by the so-called \textsc{L\'{e}vy-Khintchine} formula (see \cite{45.},  \cite{SATO}), while $({\bf Y_t})_{t\in[0,1)}$ is continuous in probability, $(\mathcal{A}_t)_{t\in[0,1)}-$ adapted, and the random variable ${\bf Y}_t-{\bf Y}_s$ is $\mathcal{P}-$ independent to $\mathcal{A}_s$ with law $({\bf Y_{t-s}}-{\bf Y_0})_\star\mathcal{P}\in {\bf M}_1(\mathbb{R}^d)$, whenever $0\leq s<t<1$. Further assume that the \textsc{L\'{e}vy} measure of $({\bf Y_t})$ satisfies 
  the integrability condition \begin{equation} \label{integrcond} \int_{\mathbb{R}^d\setminus {B}({\bf 0_{\mathbb{R}^d}}, 1)} \|{\bf x}\|^2_{\mathbb{R}^d} m({\bf dx}) <+\infty; \end{equation}  in particular, those processes 
  encompass those among compound \textsc{Poisson} processes on time interval $[0,1)$, whose \textsc{L\'{e}vy} measure satisfy~(\ref{integrcond}), notably the simple \textsc{Poisson} process. For ${\bf x} \in \mathbb{R}^d$, define $${\bf X_t} = {\bf x} + {\bf B_t} + \int_0^t{\bf Y_s} ds, \ \forall t\in [0,1], \ \mathcal{P}-a.s.,$$ and $\eta= {\bf X}_{\star}\mathcal{P}\in {\bf M}_1(W)$. Then, from example  25.12 p.163 of \cite{SATO}, we obtain $$\mathbb{E}_{\mathcal{P}}\left[\int_0^1\|{\bf Y_t} \|^2_{\mathbb{R}^d} dt\right] <+\infty,$$ and $\mathbb{E}_{\mathcal{P}}[\|{\bf Y_t} \|_{\mathbb{R}^d}] <+\infty$, $\forall t\in [0,1)$. Since the random variable ${\bf Y_t}-{\bf Y_s}$ is $\mathcal{P}-$ independent to $\mathcal{A}_s$, for any $0\leq s< t <1$, we deduce that $({\bf Z_t})_{t\in[0,1)}$ is an $(\mathcal{A}_t)_{t\in[0,1)}$ martingale, where ${\bf Z_t}(\omega)= {\bf Y_t}(\omega) -\mathbb{E}_{\mathcal{P}}[{\bf Y_t}]$, $\forall \omega \in \Omega$, $\forall t\in [0,1)$ (see Proposition 3.17, p.97 of \cite{14.}), so that Proposition~\ref{FBPROPMV} applies. As a consequence, if we define $$\phi(\nu) = \begin{cases} \mathbb{E}_\nu\left[\int_0^1\frac{\|{\bf v_s^\nu}\|^2_{\mathbb{R}^d}}{2} ds \right] \ if \  \nu \in {\bf M}_1^{\mathbb{S}}(W) \\  +\infty  \ if \ \nu \in   {\bf M}_1(W)\setminus  {\bf M}_1^{\mathbb{S}}(W) \end{cases},$$
 then $\eta={\bf X}_{\star}\mathcal{P}\in {\bf M}_1(W)$ satisfies $\delta \phi_{\eta}[{\bf h}] =0 ,$ for all ${\bf h}\in L^2_a(\eta,H^1)$ such that 
  $h_0=h_1 =0, \ \eta-a.s. ,$ and $\int_W {\bf h} \ d{\eta} ={\bf 0}_{H^1}$.
    \end{example}

    \begin{example}  In the one dimensional case $d=1$, denote by $\nu_0= \delta^{Dirac}_0\in {\bf M}_1(\mathbb{R})$ the \textsc{Dirac} mass concentrated at $0\in \mathbb{R}$, and denote by $\nu_1= \mathcal{N}\left(3- \exp(1), \frac{7}{3}\right) \in {\bf M}_1(\mathbb{R})$ the \textsc{Gauss} probability distribution on $(\mathbb{R}, \mathcal{B}_{\mathbb{R}})$ with density $$\rho_{\nu_1}(x)= \sqrt{\frac{3}{14 \pi}} \exp\left(- 3\frac{(x - 3 + \exp(1))^2}{14} \right), \ \forall x\in \mathbb{R},$$ with respect to the \textsc{Lebesgue} measure. Define   $v :  [0,1] \times W \to \mathbb{R}$ by $$v_t({\bf \omega}) = \omega(t) + \exp(-t) - \int_0^t \exp(s-t)(\omega(s) + \exp(s)) ds, $$ for all $(t,{\bf \omega})\in  [0,1] \times W.$ Using \textsc{Gr\"{o}nwall}'s lemma (see \cite{48.}), it is an easy task to check that the stochastic differential equation $$ dX_t = dB_t + v_t(X) dt; X_0=0 $$ has a unique strong solution (see \cite{45.}). Let $(\Omega,\mathcal{A},(\mathcal{A}_t)_{t\in [0,1]}, \mathcal{P})$ be a complete stochastic basis, where there exists an $(\mathcal{A}_t)-$Brownian motion $(B_t)_{t\in[0,1]}$, denote by  $(X_t)_{t\in [0,1]}$ the strong solution of this stochastic differential equation with respect to this Brownian on this complete stochastic basis, and  set $Y_t= v_t(X)$, $\forall t\in[0,1)$, $\mathcal{P}-a.s.$; $(X_t)$ depends explicitly on past values of $(B_t)$, since the \textsc{It\^{o}} calculus ensures that $X_t= B_t + \int_0^t\left( B_s + 2 - \exp(s)\right) ds$, $\forall t\in [0,1]$, $\mathcal{P}-a.s.$. Then, there exists an $(\mathcal{A}_t)-$martingale $(Z_t)_{t\in[0,1)}$ on this complete stochastic basis, such that the $(\mathcal{A}_t)-$ adapted processes $(X_t)$ and $(Y_t)$ satisfy the forward-backward system $$dX_t = dB_t + Y_t dt,$$ and $$Y_t = Z_t  + \int_{\mathbb{R}} y \  {Y_t}_\star \mathcal{P}(dy),$$ with constraints $${X_0}_\star \mathcal{P} =  \nu_0, \ and \ {X_1}_\star \mathcal{P} = \nu_1.$$ Moreover, denote by $\nu_{crit}= {\bf X}_\star \mathcal{P} \in {\bf M}_1(W)$ the unique law of solutions to this stochastic differential equation, whose entropy with respect to the \textsc{Wiener} measure is finite (see \cite{17.}). Setting $\phi$ as in Example~\ref{ex2}, we conclude similarly that $\nu_{crit}$ is extremal for $\phi$ with respect to those average preserving variations.
\end{example}

\end{document}